\documentclass[12pt]{amsart}

\usepackage{amssymb}
\newcommand{\Ceps}{\mathbb{C}(\!(\epsilon)\!)}
\newcommand{\Sp}{\mathrm{Sp}}
\newcommand{\Cdot}{{{\text{\Large$\cdot$}}}}
\newcommand{\CC}{\mathbb C}
\newcommand{\RR}{\mathbb R}
\newcommand{\ZZ}{\mathbb Z}

\newcommand\A{\mathcal A}
\newcommand\pol{\mathrm{pol}}
\newcommand{\dontprint}[1]{\relax}
\newcommand{\Curv}{C}
\newcommand{\ev}{\mathrm{ev}_1}
\newcommand{\pr}{\mathit{pr}}
\newcommand{\HH}{\mathit{HH}}

\newcommand{\ii}{\hspace{1pt}\mathrm{i}\hspace{1pt}}
\newtheorem{thm}{Theorem}[section]
\newtheorem{lemma}{Lemma}[section]
\newtheorem{proposition}{Proposition}[section]
\newtheorem{cor}{Corollary}[section]
\newenvironment{remark}{\par\medskip\noindent{\bf Remark.}}
{\par\medskip}
\begin{document}
\title[Cohomology of the Weyl algebra]
{Hochschild cohomology of the Weyl algebra and
traces in deformation quantization}
\author{Boris Feigin}
\address{Landau Institute for Theoretical Physics,
Moscow}

\email{}

\author{Giovanni Felder}
\address{Department of Mathematics,
ETH Zurich, CH-8092 Zurich}

\email{felder@math.ethz.ch}

\author{Boris Shoikhet}
\address{Department of Mathematics,
ETH Zurich, CH-8092 Zurich}

\email{borya@math.ethz.ch}

\begin{abstract}
We give a formula for a cocycle generating the Hoch\-schild
cohomology of the Weyl algebra with coefficients in its dual.
It is given by an integral over the configuration space
of ordered points on a circle. Using this formula and
a non-commutative version of formal geometry, we obtain
an explicit expression for the canonical trace in deformation
quantization of symplectic manifolds.
\end{abstract}

\maketitle
\section{Introduction}\label{sect1}
Let $(M,\omega)$ be a connected
symplectic manifold of dimension $2n$. 
In deformation quantization one considers formal local
associative deformations (star-products)
$f\star g=fg+\frac\epsilon2\{f,g\}
+\epsilon^2B_2(f,g)+\cdots$ of the product of smooth
complex-valued
functions whose first order term is given by the
Poisson bracket associated to $\omega$. Such deformations
were classified by Fedosov (see \cite{F}) up to
natural equivalence: they are in one to one correspondence
with their {\em characteristic classes}:
formal series $\Omega=\Omega_0+\epsilon\Omega_1+\cdots$
of cohomology classes $\Omega_i\in H^2(M,\mathbb C)$,
such that $-\Omega_0$ is the class of the symplectic
form. An important property of these deformations is
that they admit a {\em trace} (a linear functional
$\tau$ on compactly supported functions with values in 
$\Ceps$ such that $\tau(f\star g)
=\tau(g\star f)$), which is unique up to normalization.
Moreover there is a canonical normalization given
by a local condition: locally all deformations are
equivalent to the Weyl algebra, and on the Weyl
algebra there is a canonical trace given by the
integral with respect to the Liouville measure
(divided by $(2\pi \ii\epsilon)^{n}$). So
one requires that the trace coincides with this
canonical trace when restricted to functions with
support in the neighborhood of any point of $M$.
So to each star-product we have a unique trace 
obeying this local normalization condition. It is
called canonical trace.

Interesting formulas are then obtained
if we  evaluate the trace on 
general (compactly supported) functions on $M$: already
in the simplest case of the constant function $1$ on
a compact symplectic manifold the answer is non-trivial.
Namely, one has the following special case of the
Fedosov/Nest-Tsygan index theorem, see \cite{F},
\cite{NT1}.

\begin{thm}\label{t-1}(Fedosov, Nest--Tsygan)
Suppose $(M,\omega)$ is a compact symplectic manifold and let $\star$ be a star-product
with characteristic class $\Omega=-\omega+\epsilon\Omega_1
+\epsilon^2\Omega_2+\cdots$. Then the canonical trace
associated to $\star$ obeys
\begin{equation}\label{e-1}
\mathrm{Tr}(1)=\frac1{(2\pi\ii)^n}
\int_M\hat{\mathrm A}(\mathit{TM})
\exp\,(-\Omega/\epsilon).
\end{equation}
\end{thm}
Here $\hat{\mathrm A}(\mathit{TM})$ is a characteristic class of
the tangent bundle $\mathit{TM}$. A de Rham representative 
of this class is
$\mathrm{det}^{\frac12}\frac{ R/2}{\sinh(R/2)}$ where
$R$ is the curvature of any %symplectic
connection on $\mathit{TM}$.

The purpose of this paper 
is to generalize this theorem and give a similar
formula for $\mathrm{Tr}(f)$, as an integral
$\int_M\psi(f)$
of a closed differential $2n$-form $\psi(f)$, 
for any function $f$.

The main idea is to apply a non-commutative version of 
the ``formal differential geometry'' of \cite{GK}.
In formal differential geometry one considers, given
a $d$-manifold $M$, the
infinite dimensional manifold $\tilde M$
of formal parametrizations
($\infty$-jets at $0$ 
of local diffeomorphisms $\mathbb R^{d}
\to M$) as a principal homogeneous space of the
Lie algebra $W_d$ of formal vector fields, 
acting by infinitesimal change of coordinates.
The action of the subalgebra $\mathfrak{gl}_d$ of
linear vector fields integrates to an action of
$\mathrm{GL}_d$
and $\tilde M/\mathrm{GL}_d$ is homotopy equivalent to
$M$: there is a section $M\to \tilde M/\mathrm{GL}_d$
of the canonical projection 
$\tilde M/\mathrm{GL}_d\to M$ and 
all sections are homotopic to each other. 
We then have a map from the Lie algebra cohomology
of $W_d$ to the cohomology of $\tilde M$ and from
the relative Lie algebra cohomology of the pair
$(W_d,\mathfrak{gl}_d)$ to the cohomology of $\tilde M/
\mathrm{GL}_d$ and thus to the cohomology of $M$.

In our case, we should use a non-commutative
 version of this
construction. To do this, it is useful
to regard a  formal parametrization sending
$0$ to $x\in M$ as an algebra isomorphism 
from the algebra of  $\infty$-jets  of functions 
at a point $x\in M$ to the algebra 
$\mathbb C[[y_1,\dots,y_d]]$ of $\infty$-jets of
functions at $0\in\mathbb R^d$. This notion has
a natural non-commutative generalization. Namely,
on $\mathbb R^{2n}$ with standard symplectic structure
we have a canonical star-product, the Moyal product
(see \ref{ss-Hugh}). It is well-defined on formal
power series $\mathbb C[[y_1,\dots,y_{2n}]]$. The Weyl
algebra $\mathcal A_{2n}$
is the algebra of formal power series in $y_i $
with coefficients in $\mathbb C[[\epsilon]]$ and Moyal
product. A non-commutative version of a formal
parametrization at $x\in M$ is then an isomorphism
from the algebra of $\infty$-jets of functions on
$M$ at $x$ with star-product to the Weyl algebra 
$\mathcal A_{2n}$. The Lie algebra of vector fields
is then replaced by the Lie algebra
$\mathrm{Der}(\mathcal A_{2n})$ of derivations of
$\mathcal A_{2n}$ which is isomorphic 
to the quotient of $\mathcal A_{2n}$,
with Lie bracket $[a,b]_\epsilon=\frac1\epsilon
(a\star b-b\star a)$, by its center $Z_{2n}$, consisting of
constant functions. The subalgebra $\mathfrak{gl}_d$
is replaced by the Lie algebra $\mathfrak{sp}_{2n}$ of
the group of automorphisms of $\mathcal A_{2n}$ induced
by linear symplectic coordinate transformations.  

We are then led to study the Lie algebra cohomology
of the pair 
$(\mathrm{Der}(\mathcal A_{2n}),
\mathfrak{sp}_{2n})$, or equivalently $(\mathcal A_{2n},
\mathfrak{sp}_{2n}\oplus Z_{2n})$. Now this cohomology
(with trivial coefficients) has a canonical class,
described in \cite{FT2}: it comes from the
Hochschild cohomology of the polynomial 
Weyl algebra $\mathcal A^\mathrm{pol}_{2n}$ with
coefficients in the the dual bimodule
$(\mathcal A^\mathrm{pol}_{2n})^*$, which is one-dimensional
and generated by a canonically normalized
class $\tau_{2n}$ of degree $2n$. This class maps to a 
cohomology class $\Theta_{2n}$
of $(\mathcal A^\mathrm{pol}_{2n},
\mathfrak{sp}_{2n}\oplus Z_{2n})$
with values in $(\mathcal A^\mathrm{pol}_{2n})^*$.
We describe this construction below.
In \cite{FT2} an algebraic
version of the Riemann--Roch--Hirzebruch formula is given. 
In its simplest form it gives a formula for
 $\ev\Theta_{2n}
\in 
H^{2n}(\mathcal A_{2n},\mathfrak{sp}_{2n}\oplus Z_{2n})$,
the value of
$\Theta_{2n}$ on $1\in\mathcal 
A^\mathrm{pol}_{2n}$ in terms
of ad-invariant polynomials on $\mathfrak{sp}_{2n}\oplus Z_{2n}$ via a Chern--Weil homomorphism. After applying
the formal geometry map we obtain a class in $H^{2n}(M,
\Ceps)$, which
is precisely the class appearing in \eqref{e-1},
see \ref{ss-Clarissa} below. 

It is then natural to expect that by evaluating 
$\Theta_{2n}$ on more general elements of $\mathcal A_{2n}$
we get, via formal geometry, classes which upon integration
give the trace of general functions. This is
indeed true and the content of one of our main results
(Theorem \ref{t-2}). 

Another result of this paper is an explicit formula
for a cocycle $\tau_{2n}$
representing $[\tau_{2n}]$ and 
mapping to a relative Lie algebra cocycle $\Theta_{2n}$,
which was hitherto only known to exist by indirect
arguments. In particular, this leads in principle
to an explicit description of $\psi(f)$. Also
a direct proof of the Riemann--Roch--Hirzebruch
formula of \cite{FT2} is obtained and is presented
in Sect.~\ref{sect4}.

The explicit formula for $\tau_{2n}$, which
we present in Sect.~\ref{sect1} is given in
terms of integrals over configuration spaces of
points on a circle, and was discovered as a 
consequence of the proof by Feynman diagrams
given in \cite{Sh} of the Tsygan formality conjecture.

\section{Hochschild cohomology of the Weyl algebra}
\label{sect2}
\subsection{The Weyl algebra}\label{ss-Hugh}
The polynomial Weyl algebra $\A^\pol_{2n}$ over the
ring $K=\CC[\epsilon,\epsilon^{-1}]$ is
 the space of polynomials 
$K[p_1,\dots,p_n,q_1,\dots,q_n]$
with  the Moyal product associated with
the bivector
\[
\alpha=\frac12\sum_{i=1}^n\left(
\frac{\partial}{\partial {p_i}}\otimes\frac{\partial}{\partial {q_i}}
-
\frac{\partial}{\partial {q_i}}\otimes\frac{\partial}{\partial {p_i}}
\right)\in \mathrm{End}_\CC(\mathcal A^\pol_{2n}\otimes
 \mathcal A^\pol_{2n}).
\]
The Moyal product is given by the formula
\[
f\star g=m(\exp(\epsilon\alpha)(f\otimes g)),
\]
where $m(f\otimes g)=fg$ is the standard 
commutative product on polynomials.

This algebra is generated by $p_i,q_i$ with defining
relations $p_i\star q_j-q_j\star p_i=\epsilon
\delta_{ij}$ (thus $\hbar=\ii\epsilon$ in
the notation of physics). It is isomorphic
via $q_i\mapsto x_i$, $p_i\mapsto\epsilon
 \partial/\partial x_i$
to the algebra of differential operators in $n$ variables
with coefficients in $K[x_1,\dots,x_n]$.

The Weyl algebra is $\ZZ$-graded for the assignment
\[
\mathrm{deg}\,q_i=\mathrm{deg}\,p_i=1,\qquad
\mathrm{deg}\,\epsilon=2.
\]
We also consider the completed
Weyl algebra $\mathcal A_{2n}$
over the field of formal Laurent series $\Ceps$. As a vector space 
over $\Ceps$, 
$\mathcal A_{2n}=\CC[[p_i,q_i]](\!(\epsilon)\!)$. The
product is defined by the same Moyal formula. The completed
Weyl algebra contains the subring
$\mathcal A^+_{2n}=\CC[[p_i,q_i,\epsilon]]$, an algebra
over $\CC[[\epsilon]]$. The filtration of $\mathcal A^+_{2n}$ by degree
 defines a topology on $\mathcal A^+_{2n}$, $\mathcal A_{2n}$.

The symplectic group $\Sp_{2n}$ 
of linear transformations preserving the tensor $\alpha$ acts on
$\mathcal A^\pol_{2n},\mathcal A_{2n}^+$ and $\mathcal A_{2n}$
 by automorphisms.

We often view $\mathcal A^\pol_{2n}$ and $\mathcal A_{2n}$ as 
Lie algebras with bracket
\[
[f,g]_{\epsilon}=\frac1\epsilon(f\star g-g\star f).
\]
The map $f\mapsto [f,\ ]_\epsilon$ from $\mathcal A^\pol_{2n}$ to
the Lie algebra $\mathrm{Der}(\mathcal A^\pol_{2n})$ of derivations of
the Weyl algebra defines an exact sequence of Lie algebra homomorphisms
\[
0\to K\to \mathcal A^\pol_{2n}\to 
\mathrm{Der}(\mathcal A^\pol_{2n})\to 0.
\]
The Lie algebra $\mathfrak{sp}_{2n}$ of $\mathrm{Sp}_{2n}$
may be identified with the Lie algebra of homogeneous
quadratic polynomials in $\mathcal A^\pol_{2n}$. 
It acts on $\mathcal A^\pol_{2n}$ by the corresponding
inner derivations. 

More explicitly, introduce coordinates 
$y_i$, $i=1,\dots 2n$ by $y_{2i-1}=p_i$, $y_{2i}=q_i$
so that the standard symplectic form 
on $\mathbb R^{2n}$ is $\omega^\circ=\sum dq_i\wedge dp_i
=\frac12\sum\omega^\circ_{ij}dy_i\wedge dy_j$.
Then $\mathfrak{sp}_{2n}$
consists of $2n\times 2n$ matrices $A=(a^i_j)$
such that $a_{ij}=_{\mathrm{def}}
\sum_{k=1}^{2n}\omega^\circ_{ik}a^{k}_j$ 
is symmetric. The derivation 
of $\mathcal A^\pol_{2n}$ corresponding to $A\in\mathfrak{sp}_{2n}$ is $f\mapsto Af$ with
$Af(y)=\frac d{dt}f(\exp(-tA)y)|_{t=0}$. Then
\[
Af=[\tilde A,f]_\epsilon, \qquad
\tilde A=\frac12 \sum_{i,j=1}^{2n} a_{ij}y_iy_j.
\]

\subsection{Hochschild homology and cohomology}
The Hochschild homology $\HH_\Cdot(A)$ (with coefficients
in the bimodule $A$) of an  associative
algebra $A$ 
with unit $1$ over a commutative ring $K$ is the homology 
of the Hochschild chain complex 
$\dots\to C_2(A)\to C_1(A)\to
C_0(A)\to 0$ with $C_k(A)=A^{\otimes (k+1)}$ (the tensor product
over $K$ of $k+1$ copies of $A$) and
differential
\begin{eqnarray*}
d_{\mathrm{Hoch}} (a_0\otimes a_1\otimes\cdots\otimes a_k)
&\!=\!&
a_0\cdot a_1\otimes a_2\otimes\cdots\otimes a_k\\
&&+\sum_{j=1}^{k-1}(-1)^ja_0\otimes\cdots\otimes a_j\cdot a_{j+1}\otimes\cdots\otimes a_k
\\
&&+(-1)^ka_k\cdot a_0\otimes\cdots\otimes a_{k-1}.
\end{eqnarray*}
A more convenient complex computing this homology is
the {\em normalized chain complex} 
$\dots\to \bar C_2(A)\to \bar C_1(A)\to
\bar C_0(A)\to 0$ with $\bar C_k(A)=A\otimes 
(A/K1)^{\otimes k}$. The differential passes to
the quotient and the canonical projection $C_\Cdot(A)
\to \bar C_\Cdot(A)$ induces an isomorphism in
homology (see, e.g., \cite{M}, Chapter X, Cor.\ 2.2)
% or \cite{L}, Proposition 1.6.5).

\begin{thm}
{\cite{FT1}}
The Hochschild homology of $\mathcal A^\pol_{2n}$
is zero except in degree $2n$; $\HH_{2n}(\mathcal A^\pol_{2n})
\simeq\CC[\epsilon,\epsilon^{-1}]$
is generated by the class of the cycle
$
c_{2n}=
\sum_{\sigma\in S_{2n}}\mathrm{sign}(\sigma)\,
1\otimes y_{\sigma(1)}\otimes \cdots\otimes y_{\sigma(2n)}
\in \bar C_{2n}(\mathcal A^\pol_{2n})$, 
$y_{2i-1}=p_i$, $y_{2i}=q_i$.
\end{thm}

The Hochschild cohomology $\HH^{\Cdot}(A)=\HH^{\Cdot}(A;A^*)$
of an algebra is defined as
the cohomology of the dual complex $\mathrm{Hom}_K(\bar C(A),K)$ 
(it is the Hochschild cohomology of the algebra $A$ 
with coefficients in the bimodule $A^*$, rather than in the bimodule $A$
appearing in deformation theory). Thus $\HH^\Cdot(\mathcal
A^\pol_{2n})$ is one-dimensional and concentrated in degree
$2n$. Our aim is to write a formula for a non-trivial
$2n$-cocycle. 
\subsection{Formula for the Hochschild cocycle}
\label{sus-6}
We first introduce some notation: let $\Delta_k\subset\RR^k$
be the simplex $0\leq u_1\leq\cdots\leq u_k\leq1$ with
the standard orientation. For $0\leq i\neq j\leq 2n$ denote by
$\alpha_{ij}=-\alpha_{ji}\in \mathrm{End}
((\mathcal A^\pol_{2n})^{\otimes(k+1)})$ the map
\begin{eqnarray*}
\alpha_{ij}(a_0\otimes\cdots\otimes a_k)&=&
\frac12\sum_{r}\,
(
a_0\otimes\cdots\otimes
\partial_{p_{r}}a_i\otimes
\cdots\otimes
\partial_{q_{r}}a_j\otimes
\cdots\otimes a_k\\
&&-\,
a_0\otimes\cdots\otimes
\partial_{q_{r}}a_i\otimes\cdots\otimes
\partial_{p_{r}}a_j\otimes\cdots\otimes a_k)
\end{eqnarray*}
Define $\pi_{2n}\in \mathrm{End}(A^{\otimes{2n+1}})$ by
\[
\pi_{2n}(a_0\otimes\cdots\otimes a_{2n})
=\sum_{\sigma\in S_{2n}}
\mathrm{sign}(\sigma)
a_0\otimes\partial_{y_{\sigma(1)}}a_1
\otimes\cdots\otimes
\partial_{y_{\sigma(2n)}}a_{2n}.
\]
Here $y_{2i-1}=p_i$ and $y_{2i}=q_i$, $1\leq i\leq n$.
Finally, let $\mu_k:(\mathcal A^\pol_{2n})^{\otimes{(k+1)}}\to \CC$ denote the
map $a_0\otimes\cdots\otimes a_k\mapsto a_0(0)\cdots a_k(0)$,
where $a(0)$ is the constant term of the power
series $a\in \mathcal A^\pol_{2n}$.
Both $\alpha_{ij}$ and $\pi_{2n}$ commute with the action of
$\mathfrak{sp}_{2n}(\CC)$,
and $\mu_k$ is an $\mathfrak{sp}_{2n}(\CC)$-invariant
linear form.
Moreover $\pi_{2n}$ induces a well-defined map
$\bar C_{2n}(\mathcal A^\pol_{2n})=\mathcal A^\pol_{2n}\otimes(\mathcal A^\pol_{2n}/\CC1)^{\otimes{2n}}\to
(\mathcal A^\pol_{2n})^{\otimes(2n+1)}$.

The formula for a cochain which, as we will see,
is a cocycle, is
\begin{equation}\label{e-Hoco}
\tau_{2n}(a)=\mu_{2n}\int_{\Delta_{2n}}
\prod_{0\leq i<j\leq 2n}e^
{\epsilon(2u_i-2u_j+1)\,\alpha_{ij}}\pi_{2n}(a)\,
du_1\wedge\cdots\wedge du_{2n},
\end{equation}
with the understanding that $u_0=0$. Here the exponential function
is to be expanded in a Taylor series (only finitely many terms
contribute non-trivially) and integrated term by term. Note that
the maps $\alpha_{ij}$ commute so there is no ordering ambiguity.

\begin{remark} The integrals in the formula
for $\tau_{2n}$ may be ``computed''. For example,
if $n=1$ we have:
\[
\tau_2=\mu_{2}\circ F(\epsilon\alpha_{01},\epsilon\alpha_{12},
\epsilon\alpha_{20})
\circ\pi_{2},
\]
where the formal power series $F\in\CC[[x,y,z]]$ is given
by
\[
F(x,y,z)=-\frac{
(x-y)\,
e^{x+y-z}
+
(y-z)\,
e^{y+z-x}
+
(z-x)\,
e^{z+x-y}
}
{4 (x-y)\,
 (y-z)\,(z-x)
}.
\]
\end{remark}
\begin{thm}\label{t-sp}\ 

\begin{enumerate}
\item[(i)] The Hochschild cochain
$\tau_{2n}$ is an $\mathfrak{sp}_{2n}(\CC)$-invariant
cocycle in the normalized Hochschild complex
$\bar C^{2n}(\mathcal A^\pol_{2n})$, i.e., a cocycle so that
\[
\sum_{i=0}^{2n}\tau_{2n}(a_0\otimes\cdots\otimes
[a,a_i]_\epsilon \otimes\cdots\otimes a_{2n})=0,
\]
for all $a\in \mathfrak{sp}_{2n}(\CC)\subset \mathcal A^\pol_{2n}$.
\item[(ii)] $\tau_{2n}(c_{2n})=1$.
\item[(iii)] If $a\in \mathfrak{sp}_{2n}(\CC)$, then
\[
\sum_{i=1}^{2n}(-1)^{i}
\tau_{2n}(a_0\otimes a_1\otimes\cdots\otimes
a_{i-1}\otimes a\otimes a_{i}\otimes
\cdots\otimes a_{2n-1})=0.
\]
\end{enumerate}
\end{thm}
\medskip

This theorem is proved in \ref{ss-jupiter} below.

\begin{remark} {\it Comparison with Tsygan formality.}
Theorem \ref{t-sp} (i) is a consequence of the
formality theorem for Hochschild chains 
(Tsygan conjecture), proved in \cite{Sh} in terms of
integrals over configuration spaces, and the
the formula for $\tau_{2n}$ was found 
by simplifying these integrals. 
Here is how the story goes. First of all we have
Kontsevich's formality
theorem \cite{K}, which gives an $L_\infty$-%
quasi\-iso\-mor\-phism between the differential graded
Lie algebra of (formal) polyvector fields on $\RR^d$
(with trivial differential)
and the differential graded Lie algebra of
continuous Hochschild
cochains of $\mathcal O_d=\CC[[x_1,\dots,x_d]]$ with
values in $\mathcal O_d$ and with the Gerstenhaber bracket.
Now differential forms on $\RR^d$ form a module over
polyvector fields (acting by Lie derivative) and
Hochschild chains form a module over Hochschild cochains.
Both modules can be considered as $L_\infty$-modules
over polyvector fields via Kontsevich's formality.
Based on this observation, Tsygan \cite{T} conjectured
the existence of a morphism of $L_\infty$-modules
from the module of Hochschild chains of $\mathcal O_d$
to the module
of differential forms. This morphism was constructed
in terms of integrals over configuration spaces in
\cite{Sh}, and Tsygan's conjecture was proved there using
Stokes' theorem.

This morphism of $L_\infty$-modules may be linearized
at any solution of the Maurer--Cartan equation in
the Lie algebra of polyvector fields, yielding
a tangent map between cohomology spaces. In particular,
the bivector $\alpha$ is a Poisson bivector field
and thus a solution of the Maurer--Cartan
equation for $d=2n$. The corresponding tangent map is an
isomorphism from the Hochschild homology of $\mathcal A^\pol_{2n}$ to
the homology of the complex of differential forms with
differential $L_\alpha$, the Lie derivative in the
direction of the bivector field $\alpha$. The latter
homology is zero except in degree $2n$, where it is
spanned by the class of $dy_1\wedge\cdots\wedge dy_{2n}$.
The dual map to the tangent map
sends the dual cocycle to $\tau_{2n}$. Its coefficients
are given by integrals over the configuration space
of $2n+1$ distinct points on the boundary of a disk and an
arbitrary number of distinct points in its interior.
The integral over the interior points may be performed
explicitly and one obtains the formula \eqref{e-Hoco}.

On the other hand, the important
statement (iii) comes somewhat
as a surprise in this approach.
\end{remark}

\subsection{Proof of Theorem \ref{t-sp}}\label{ss-jupiter}

\begin{proposition}\label{l-101}
 $\tau_{2n}$ is a cocycle in $\bar C(\mathcal A^\pol_{2n})$
and $\tau_{2n}(c_{2n})=1$.
\end{proposition}

\medskip
\noindent{\it Proof:} Let us first check the
normalization condition $\tau_{2n}(c_{2n})=1$.
We have $\pi_{2n}(c_{2n})=(2n)!1\otimes\cdots\otimes 1$, so
$\alpha_{ij}$ acts trivially and
the integral reduces to the volume $1/(2n)!$ of $\Delta_{2n}$. Finally, $\mu_{2n}(1\otimes\cdots\otimes 1)=1$.

The proof of the cocycle property
is based on Stokes' theorem.
We introduce a closed $2n$-form $\eta$ on $\RR^{2n+1}$
with values in $((\mathcal A^\pol_{2n})^{\otimes(2n+2)})^*$.
The integral of $0=d\eta$
on the simplex
$\Delta_{2n+1}$, calculated by Stokes' theorem, yields
the cocycle condition. The construction of $\eta$ goes as
follows. Let $\pi_{2n}^i\in\mathrm{End}((\mathcal A^\pol_{2n})^{\otimes(2n+2)})$ be
the map
\[
\pi_{2n}^i
=\sum_{\sigma\in S_{2n}}
\mathrm{sign}(\sigma)
\mathrm{Id}\otimes\partial_{y_{\sigma(1)}}\otimes\cdots
\otimes\partial_{y_{\sigma(i-1)}}
\otimes \mathrm{Id}
\otimes\partial_{y_{\sigma(i)}}
\otimes\cdots
\otimes\partial_{y_{\sigma(2n)}}
\]
Then we set (the caret denotes omission)
\[
\eta=\mu_{2n+1}\sum_{i=1}^{2n+1}
\prod_{0\leq j<k\leq 2n+1}e^{
(2u_j-2u_k+1)\,\alpha_{jk}
}
du_1\wedge\cdots\hat{du_i}\cdots\wedge du_{2n+1}
\,\pi_{2n}^i,
\]
where we set $u_0=0$.
We first show that $\eta$ is a closed form. We have
\begin{eqnarray*}
d\eta&=&\mu_{2n+1}\sum_{i=1}^{2n+1}
(-1)^{i}\sum_{l=0}^{2n+1}2\,\alpha_{li}\,\pi_{2n}^i\\
&&
\prod_{0\leq j<k\leq 2n+1}e^{
(2u_j-2u_k+1)\,\alpha_{jk}
}
du_1\wedge\cdots\wedge du_{2n+1}.
\end{eqnarray*}
We claim that $\phi\equiv\mu_{2n+1}\sum_{i=1}^{2n+1}
(-1)^i\sum_{l=0}^{2n+1}2\,\alpha_{li}\,\pi_{2n}^i=0$, from which
it follows that $d\eta=0$.
To prove the claim, notice that
\[
\phi(a_0\otimes\cdots\otimes a_{2n+1})=\sum_{r=1}^n
(\partial_{y_{2r-1}}b_{2r}-\partial_{y_{2r}}b_{2r-1})(0),
\]
where we set
\[
b_j=\sum_{i=1}^{2n+1}(-1)^i\sum_{\sigma\in S_{2n}}
a_0
\frac{\partial a_1}{\partial y_{\sigma(1)}}
\cdots
\frac{\partial a_{i-1}}{\partial y_{\sigma(i-1)}}
\frac{\partial a_{i}}{\partial y_{j}}
\frac{\partial a_{i+1}}{\partial y_{\sigma(i)}}
\cdots
\frac{\partial a_{2n+1}}{\partial y_{\sigma(2n)}}\,.
\]
This expression is $a_0$
times the Lagrange
expansion with the respect to the first column
$(\partial a_i/\partial y_j)_{i=1,
\dots,2n+1}$
of the determinant of a $2n+1$ by $2n+1$ matrix
whose remaining columns are $(\partial a_i/\partial y_k)$,
$k=1,\dots,2n$. Since two columns are equal, the
determinant vanishes.
 Thus $\eta$ is closed.

We now apply Stokes' theorem
\[
0=\sum_{k=0}^{2n+1}(-1)^k\int_{\Delta_{2n}}\iota_k^*\eta.
\]
Here $\iota_k:\Delta_{2n}\to\Delta_{2n+1}$ is the
$k$th face map:
\begin{eqnarray*}
\iota_0(u_1,\dots,u_{2n})&=&(0,u_1,\dots,u_{2n}),\\
\iota_k(u_1,\dots,u_{2n})&=&(u_1,\dots,u_k,u_k,\dots,u_{2n}),
\qquad 1\leq k\leq 2n,\\
\iota_{2n+1}(u_1,\dots,u_{2n})&=&(u_1,\dots,u_{2n},1).
\end{eqnarray*}
The pull-backs  $\iota_k^*\eta$ can then be computed.
If we evaluate them on $a=a_0\otimes\cdots\otimes a_{2n+1}$,
the result is
\[
\int_{\Delta_{2n}}\iota_k^*\eta(a)
=\left\{\begin{array}{ll}
\tau_{2n}(a_0\otimes\cdots\otimes
a_k\star a_{k+1}\otimes\cdots\otimes a_{2n+1}),&
0\leq k\leq 2n,\\
\tau_{2n}(a_{2n+1}\star a_0\otimes a_1\otimes\dots\otimes
a_{2n}),&
k=2n+1.
\end{array}\right.
\]
These are precisely the terms appearing in the
Hochschild differential. The proof is complete.
\hfill $\square$

\begin{lemma}\label{l-102}
 $\tau_{2n}$ is $\mathfrak{sp}_{2n}(\CC)$-invariant.
\end{lemma}

\noindent{\it Proof:} The claim follows immediately
from the fact that $\pi_{2n}$, $\alpha_{ij}$, and $\mu_{2n}$
are homomorphisms of $\mathfrak{sp}_{2n}(\CC)$-modules (with
trivial action of $\mathfrak{sp}_{2n}(\CC)$ on $\CC$).
\hfill $\square$

\medskip

Next, we study the behaviour of $\tau_{2n}$ under permutation
of its last $2n$ arguments. Let $\psi:\RR\to[-1,1]$ be
the function so that $\psi(u)=2u+1$ if $-1\leq u<0$
and $\psi(u+1)=\psi(u)$. Clearly, $\psi(-u)=-\psi(u)$
(for almost all $u$).
Let $\omega_{2n}$ be the piecewise smooth
$2n$-form on $\RR^n$ with values
in $\mathrm{End}((\mathcal A^\pol_{2n})^{\otimes (2n+1)})$:
\begin{equation}\label{e-omega}
\omega_{2n}=\prod_{0\leq i<j\leq 2n}\exp\left(
\psi(u_i-u_j)\,\alpha_{ij}
\right)\,|_{u_0=0}\,
du_1\wedge\cdots\wedge du_{2n}.
\end{equation}
Let the symmetric group $S_{2n}$ act on $\RR^{2n}$ by
permutation of coordinates and introduce
\[
\tau_{2n}^\sigma=\mu_{2n}\circ\int_{\sigma(\Delta_{2n})}\omega_{2n}\circ
\pi_{2n}.
\]
As $\psi(u_i-u_j)=2u_i-2u_j+1$ if $0<u_i<u_j<1$, we
have $\tau_{2n}^{\mathrm{Id}}=\tau_{2n}$.

\begin{lemma}\label{l-103}
For any permutation $\sigma\in S_{2n}$,
\[
\tau_{2n}(a_0\otimes a_{\sigma^{-1}(1)}\otimes\cdots\otimes
a_{\sigma^{-1}(2n)})=
\mathrm{sign}(\sigma)\tau_{2n}^\sigma(a_0\otimes a_1\otimes\cdots
\otimes a_{2n}).
\]
\end{lemma}

\medskip
\noindent{\it Proof:} Let $S_{2n}$ act on $(\mathcal A^\pol_{2n})^{\otimes(2n+1)}$
by permutations of the last $2n$ factors. This representation
will be denoted by
$\sigma\mapsto\rho(\sigma)\in\mathrm{End}(
(\mathcal A^\pol_{2n})^{\otimes(2n+1)})$. Clearly
\[
\pi_{2n}\circ\rho(\sigma)=\mathrm{sign}(\sigma)\rho(\sigma)\circ\pi_{2n},
\qquad \mu_{2n}\circ\rho(\sigma)=\mu_{2n}.
\]
We claim that
\[
\omega_{2n}\circ\rho(\sigma)=\mathrm{sign}(\sigma)
\rho(\sigma)\circ\sigma^*\omega_{2n}.
\]
Indeed, conjugating $\omega_{2n}$ by $\rho(\sigma)$ gives
a product ($\sigma(0)=0$)
\begin{eqnarray*}
\lefteqn{\prod_{0\leq i<j\leq 2n}\exp\left(
\psi(u_i-u_j)\,\alpha_{\sigma(i)\sigma(j)}
\right)}\\
&&=\prod_{0\leq \sigma^{-1}(i)<\sigma^{-1}(j)\leq 2n}
\exp\left(
\psi(u_{\sigma^{-1}(i)}-u_{\sigma^{-1}(j)})\,\alpha_{ij}
\right)
\\
&&=
\prod_{0\leq i<j\leq 2n}
\exp\left(
\psi(u_{\sigma^{-1}(i)}-u_{\sigma^{-1}(j)})\,\alpha_{ij}
\right).
\end{eqnarray*}
In the last step, we use the fact that $\alpha_{ij}=-\alpha_{ji}$ and $\psi(-u)=-\psi(u)$. Finally, we change variables
in the integral, taking into account an additional
$\mathrm{sign}(\sigma)$ coming from the orientation reversal.
\hfill $\square$

\medskip

We can now complete the proof of Theorem \ref{t-sp}.
The statements (i) and (ii) follow from
Proposition \ref{l-101} and Lemma \ref{l-102}.
It remains to prove (iii). By
Lemma \ref{l-103}, the left-hand side may be written as
$-
\mu_{2n}\circ\int_D\omega_{2n}\circ\pi_{2n}(a_0\otimes a\otimes a_1\otimes
\cdots\otimes a_{2n-1})$.
 The integral is over $0\leq
 u_1\leq 1$, $0\leq u_2\leq
\cdots\leq u_{2n}\leq 1$. Since $a$ is a quadratic polynomial,
only terms linear in $\alpha_{1i}$ contribute non-trivially,
and we may replace the product of
exponentials
$\prod_{i\neq 1}\exp(\psi(u_1-u_i)\alpha_{1i})$
involving $u_1$
by $\sum_{i\neq1}\psi(u_1-u_i)\alpha_{1i}$. Then the integral
over $u_1$ gives zero since
\[
\int_{0}^1\psi(u_1-u_i)du_1=0.
\]

\section{Lie algebra cohomology}
For the application to deformation quantization it is useful
to reformulate our results in terms of Lie algebra cohomology.

\subsection{Hochschild and Lie algebra cohomology}\label{ss-HLAC}
Let $K$ be a commutative ring and 
$\mathfrak g\supset \mathfrak h$ be a pair of
Lie algebras over $K$.
Let $C^\Cdot(\mathfrak g,\mathfrak h;M)
=\mathrm{Hom}(\wedge^\Cdot(\mathfrak g/\mathfrak h),M)^{\mathfrak h}$ denote
the relative Lie algebra co\-chain complex of the 
pair $(\mathfrak g,\mathfrak h)$
with coefficients in a
$\mathfrak g$-module $M$. It is the subcomplex of the
Lie algebra cochain complex $C^\Cdot(\mathfrak g;
M)=\mathrm{Hom}_{\CC}(\wedge^\Cdot\mathfrak g, M)$ consisting
of $\mathfrak h$-invariant 
cochains vanishing if any argument is in $\mathfrak h$. The
differential is given by the formula 
\begin{eqnarray*}
d_\mathrm{Lie}c(a_1\wedge\dots\wedge a_{p+1})&=&
\sum_{i=1}^{p+1}
(-1)^{i-1}a_i\cdot c(a_1\wedge\cdots\hat a_i\cdots
\wedge a_{p+1})\\
&+&\sum_{i<j}
(-1)^{i+j}c([a_i,a_j]\wedge\cdots
 \hat a_i\cdots
 \hat a_j\cdots
\wedge a_{p+1}).\\
\end{eqnarray*}
If $M$ is the trivial module $K$, then $
C^\Cdot(\mathfrak g,\mathfrak h)=
C^\Cdot(\mathfrak g,\mathfrak h;K)$
is a differential graded commutative $K$-algebra
with respect to the cup-product
\begin{eqnarray*}
c\cup c'(a_1\wedge\cdots\wedge a_{p+q})
&=&\sum\mathrm{sign}(\sigma)\,
c(a_{\sigma(1)}\wedge\cdots\wedge a_{\sigma(p)})\\
&&\cdot
c'(a_{\sigma(p+1)}\wedge\cdots\wedge a_{\sigma(p+q)}),
\end{eqnarray*}
where the sum is over $(p,q)$-shuffles,
i.e., permutations such that 
$\sigma(1)<\cdots<\sigma(p)$ and
 $\sigma(p+1)<\cdots<
\sigma(p+q)$.

For any algebra $A$ over $\CC$ with unit,
denote by  $\mathfrak{gl}_N(A)$
the Lie algebra of $N\times N$ matrices with coefficients
in $A$. It is the tensor product
$M_N(\CC)\otimes A$ of the algebra of $N\times N$ matrices
with $A$, considered as a Lie algebra. It contains the 
Lie subalgebra $\mathfrak{gl}_N=M_N(\CC)\otimes 1$.
 Then we have  chain maps $\phi^N: C^\Cdot(A)\to 
C^\Cdot(\mathfrak{gl}_N(A);
\mathfrak{gl}_N(A)^*)$:
\begin{eqnarray*}
\lefteqn{
\phi^N(\tau)
(M_1\otimes a_1,\cdots, M_k\otimes a_k)
(M_0\otimes a_0)}\\
&
=&\sum_{\sigma\in S_k}
\mathrm{sign}(\sigma)
\tau(a_0\otimes a_{\sigma(1)}\dots\otimes a_{\sigma(k)})
\mathrm{tr}(M_0M_{\sigma(1)}\cdots M_{\sigma(k)}).
\end{eqnarray*}
These maps are compatible with the inclusion 
$i_{N'N}:\mathfrak{gl}_N(A)\to \mathfrak{gl}_{N'}(A)$,
$N<N'$ obtained by embedding an $N\times N$ matrix in
the first rows and columns of an $N'\times N'$ matrix
and completing with zeros. Namely $i_{N'N}$ induces
a restriction map $i_{N'N}^*$ on complexes and
 one has
\[
\phi_{N}=i_{N'N}^*\circ \phi_{N'}.
\]
Let $S(V)=\oplus S^j(V)$ denote the symmetric algebra
of a $K$-module $V$.
By composing $\phi^N$ with the dual of the 
homomorphism of $\mathfrak{gl}_N(A)$-modules
 $S(\mathfrak{gl}_N(A))\to\mathfrak{gl}_N(A)$
\[
x_1\cdots x_k\mapsto\frac1{k!}
\sum_{\sigma\in S_k}x_{\sigma(1)}\cdots x_{\sigma(k)}
\]
(the product on the left is the product in the symmetric
algebra of the vector space $S(\mathfrak{gl}_N(A))$; the
product on the right is the associative product of 
$M_N(\CC)\otimes A$), we obtain an extension
of $\phi^N$ to 
a chain map $\phi^N_j:C^\Cdot(A)\to 
C^\Cdot(\mathfrak{gl}_N(A);
S^j(\mathfrak{gl}_N(A))^*)$, for all $j\geq1$.
These maps induce maps
\[
\phi_{k,j}^N:\HH^k(A)\to H^k(\mathfrak{gl}_N(A),S^j\mathfrak {gl}_N(A)^*)
\]
($j\geq1$) on cohomology.
\begin{thm}\label{t-Dido}
 \cite{FT2} Let $N>>n$. If $A=\mathcal A^\pol_{2n}$, then 
 $\phi^N_{k,j}$ is an isomorphism for all $0\leq k\leq 2n$
and $j\geq1$. In particular, 
\[
H^{k}(\mathfrak{gl}_N(\A_{2n}^\pol);S^j\mathfrak{gl}_N(\A_{2n}^\pol)^*)
=\left\{\begin{array}{rl}
0,&k<2n,\\
\CC[\epsilon,\epsilon^{-1}],&k=2n.
  \end{array}
\right.
\]
\end{thm}

\subsection{The Lie algebra cocycle $\Theta^N_{2n}$}
\label{ss-Juno}
Using the map $\phi^{N}_{2n,1}$ of \ref{ss-HLAC} we get 
 Lie algebra cocycles $\Theta^N_{2n}=
\phi^{N}_{2n,1}(\tau_{2n})\in
C^{2n}(\mathfrak{gl}_N(\mathcal{A});
\mathfrak{gl}_N(\mathcal{A})^*)$ compatible with 
the inclusions $\mathfrak{gl}_N\subset \mathfrak{gl}_{N'}$,
$N<N'$.

\begin{cor}\label{c-Helena}
The cocycle $\Theta^N_{2n}$ belongs to
the relative complex of the pair  
$(\mathfrak{gl}_N(\mathcal{A}^\pol_{2n}),\mathfrak{gl}_N
\oplus \mathfrak{sp}_{2n})$ and defines a non-trivial class
\[
[\Theta^N_{2n}]\in H^{2n}(\mathfrak{gl}_N(\mathcal{A}^\pol_{2n}),
\mathfrak{gl}_N
\oplus \mathfrak{sp}_{2n};\mathfrak{gl}_N(\mathcal{A}^\pol_{2n})^*),
\]
compatible with inclusions $i_{N,N'}$.
It obeys
\begin{equation}\label{e-priamus}
\Theta^N_{2n}(p_1\wedge q_1\wedge\cdots\wedge p_n\wedge q_n)(1)=N.
\end{equation}
\end{cor}

\noindent{\it Proof:}
Since $\tau_{2n}$ belongs to the normalized 
Hochschild cochain complex, its image $\Theta^N_{2n}$
lies in the relative complex $C^{2n}(\mathfrak{gl}_N(\mathcal{A}),\mathfrak{gl}_N;
\mathfrak{gl}_N(\mathcal{A})^*)$.  Theorem \ref{t-sp} (i) implies that $\Theta^N_{2n}$
is $\mathfrak{sp}_{2n}$-
invariant, (iii) that it vanishes if any argument belongs
to $\mathfrak{sp}_{2n}$; (ii) shows 
\eqref{e-priamus} and in particular that it is non-trivial
since it does not vanish on the cycle 
$p_1\wedge q_1\wedge\cdots\wedge p_n\wedge q_n\otimes 1$ 
in the relative Lie algebra chain 
complex $C_{2n}(\mathfrak{gl}_N(\mathcal A^\pol_{2n}),
\mathfrak{gl}_N\oplus \mathfrak{sp}_{2n};
\mathfrak{gl}_N(\mathcal A^\pol_{2n}))=(\wedge^\Cdot(\mathfrak g/\mathfrak h))_\mathfrak h$ dual to the cochain complex.

\section{Traces in  deformation quantization of symplectic manifolds}\label{sect3}
In this section we give a formula for the trace in
deformation quantization in terms of the Hochschild
cocycle of the previous section. We start by reviewing 
Fedosov's deformation quantization of symplectic 
manifolds, following \cite{F}.

\subsection{Star-products and quantum algebras}
\label{ss-Sally}
Let $(M,\omega)$ be a symplectic manifold. A {\em star-product} \cite{BFFLS} on $M$
is a formal associative deformation
\[
f\star g=fg+\epsilon B_1(f,g)+\epsilon^2B_2(f,g)+\cdots, f,g\in C^\infty(M),
\]
of the product of smooth functions on $M$,
so that $B_i$ are  bidifferential operators and $f\star 1=1\star f=1$.
Such deformations were classified up to a natural equivalence by Fedosov.
They are in one-to-one correspondence with series $\Omega=
\omega+\epsilon\omega_1+\epsilon^2\omega_2+
\cdots\in H^2(M,\CC{})[[\epsilon]]$ starting from
the class of the symplectic form.

Let us recall Fedosov's construction (see \cite{F}).
A {\em symplectic frame} at $x\in M$ is an ordered basis $e_1,\dots,e_{2n}$
of the tangent bundle, so that 
$\omega(e_{2i},e_{2i-1})=-\omega(e_{2i-1},e_{2i})=1$, $1\leq i\leq n$
and
$\omega(e_k,e_l)=0$ for all remaining pairs $(k,l)$.
Symplectic frames form a principal $\mathrm{Sp}_{2n}$-bundle 
$F_{\mathrm{Sp}}(M)$ over $M$.
Then $E(M)=F_{\Sp}(M)\times_{\Sp_{2n}}
\mathcal A^+_{2n}$ is a bundle 
of filtered
$\CC[[\epsilon]]$-algebras over $M$. The fiber at $x$ of this bundle
may be thought of as the quantization of the algebra of functions on
an infinitesimal neighbourhood of $0$ of the symplectic vector space $T_xM$. In particular there is a trivial
line subbundle $Z(M)$ whose fiber at $x$ consists
of  the central (constant)
functions on $T_xM$, 
and  a canonical embedding
 $j_1:T^*M\hookrightarrow E(M)$ sending
a cotangent vector to the corresponding linear function
on $\mathit{TM}$. 
The choice of a local 
trivialization $F_\Sp(M)\to M\times \Sp_{2n} $ 
identifies locally
$E(M)$ with $M\times \mathcal A^+_{2n}$. 
Sections of $E(M)$ are then locally of the
form $\sum_{I\in\ZZ_{\geq0}^{2n},k\in\ZZ_{\geq 0}} 
a_{I,k}(x)y^I\epsilon^k$, in multiindex notation.
We call a section smooth if all coefficients $a_{I,k}(x)$ are smooth functions,
and denote by $\Gamma(M,E(M))$ the algebra of smooth sections of $E(M)$. More generally, we have the algebra 
$\Omega^\Cdot(M,E(M))=\Omega^\Cdot(M)\otimes_{C^\infty(M)}
\Gamma(M,E(M))$ of smooth differential forms
with values in $E(M)$. 

The composition $\mathit{TM}\to T^*M\to E(M)$ of the isomorphism
induced by $\omega$ with the canonical embedding  $j_1$
gives
an element $A_0\in \Omega^1(M,E(M))$.

A {\em symplectic connection} on $M$ is a torsion free connection  $\nabla$ 
on $\mathit{TM}$ such that $d\omega(\xi,\eta)=\omega(\nabla\xi,\eta)+\omega(\xi,\nabla\eta)$
for any vector fields $\xi,\eta$. Symplectic connections exist on any
symplectic manifolds, see, e.g., \cite{F}. 
A symplectic connection induces a connection, also
denoted by $\nabla$, on the associated bundle $E(M)$. It can be uniquely
extended to a derivation of degree 1 
of the algebra $\Omega^\Cdot(M,E(M))$. If
 $R=\nabla^2\in
\Omega^2(M,\mathrm{End}(\mathit{TM}))$ denotes the curvature of
$\nabla$, then  
the curvature of the induced connection
$\nabla$ on $E(M)$ is $\nabla^2=[\tilde R,
\ ]_\epsilon$, where $\tilde R\in \Omega^2(M,E(M))$ is
obtained from $R$ via the embedding $\mathrm{sp}_{2n}
\to \mathcal A^+_{2n}$.

The deformation of the algebra of functions is the
algebra of horizontal sections of the bundle of algebras 
$E(M)$, with respect to a flat connection
$D:\Gamma(M,E(M))\to \Omega^1(M,E(M))$ 
of the form
 $D=\nabla+[ A,\ ]_\epsilon$ for some symplectic connection
$\nabla$ and some 1-form $A\in\Omega^1(M,E(M))$
such that $A=A_0+$ terms of degree $\geq 2$.
Note that the curvature of a connection of this form is
$D^2=[\Omega,\ ]_\epsilon$, where $\Omega\in\Omega^2(M,E(M))$
is the 
 {\em Weyl curvature} 
\[\Omega=\tilde R+
\nabla A+\frac12[A,A]_{\epsilon}.\]
The connection $D$ is flat if and only if $\Omega$
a 2-form with values in the center $Z(M)$.

\begin{thm} (Fedosov) 
\begin{enumerate}
\item For any symplectic connection and any $\Omega_\epsilon
=-\omega+\epsilon\omega_1+\dots\in\Omega^2(M,\CC{}[[\epsilon]])$ with
$d\omega_i=0$, there exists a 1-form $A\in
\Omega^1(M,\mathcal A^+_{2n})$
so that $D=\nabla +[A,\ ]$ is a flat connection with
Weyl curvature $\Omega_\epsilon$.
\item For any such flat connection on $E(M)$,  the 
algebra $\mathcal A_D(M)=\mathrm{Ker}(D)$ of horizontal section of $\Gamma(M,E(M))$
is isomorphic as a $\CC{}[[\epsilon]]$-module to $C^\infty(M)[[\epsilon]]$ and the induced product on 
$C^\infty(M)[[\epsilon]]$ is a star-product on $M$. 
\end{enumerate}
\end{thm}
Two Fedosov connections $D$, $D'$
are called {\em gauge equivalent} it there is a one-parameter family
$D_t=d+[A_t,\ ]_\epsilon$, $0\leq t\leq 1$ of Fedosov connections
interpolating between $D$ and $D'$, such that 
\[
\frac d{dt} A_t=d\lambda_t+[A_t,\lambda_t]_\epsilon,
\]
for some smooth family of sections  $\lambda_t$ of $E(M)$.
Gauge equivalent connections lead to equivalent star
products and have the same Weyl curvature.

\medskip
\noindent{\it Basic example:} Let  $M\subset\RR^{2n}$ be an open subset
with standard symplectic structure 
$\omega^\circ=\sum_{i=1}^n dx_{2i-1}\wedge dx_{2i}$. Then
$A_0=\sum\omega^\circ_{ij}y_idx_j$ and
\[
D^\circ=d+
[\sum\omega^\circ_{ij}y_idx_j,\ ]_\epsilon=\sum_{i=1}^{2n}dx_i
\left(\frac{\partial}{\partial x_i}-\frac{\partial}{\partial y_i}\right)
\]
is a Fedosov connection. Its curvature is 
$\Omega^\circ=-\omega^\circ$.
Moreover, every deformation is locally of this form, 
namely we have the ``quantum Darboux theorem'', see
\cite{F}, Theorem 5.5.1:
\begin{proposition}\label{p-0}
Let $D$ be  a Fedosov connection on $E(M)$. Then each $x\in M$ has
 an open neighbourhood $U$ such  $E(U)$, as a bundle of algebras with
connection, is isomorphic to $(E(V),D^\circ)$ for some open subset $V\subset \RR^{2n}$.
\end{proposition}

Denote by $\mathcal A_D^c(M)$ the subalgebra of compactly supported functions
in $\mathcal A_D(M)$.

\medskip

\noindent{\bf Definition.} A trace on $\mathcal A_D^c (M)$ is a
 $\Ceps$-linear map
$\mathrm{tr}:\mathcal A_D^c(M)[\epsilon^{-1}]\to \Ceps$ such
that $\mathrm{tr}(f\star g)=\mathrm{tr}(g\star f)$.

\noindent{\bf Definition.} A trace $T$ is called normalized if,
for all $f\in E_D$ with support in a sufficiently
small neighbourhood of any point,
\[
T(f)=\frac1{(2\pi\ii \epsilon)^n}
\int_{\RR^{2n}} \phi(f)\,dq_1\wedge dp_1\wedge\cdots
\wedge dq_n\wedge dp_n,
\]for any choice of a trivialization isomorphism $\phi$.

\begin{thm} Every quantum algebra $E_D(M)$ has a unique normalized trace $\mathrm{Tr}_D$.
\end{thm}

The uniqueness of such a trace 
is clear, as every function with compact support may be 
written as a finite sum of functions with arbitrarily small
support. The existence is shown in \cite{F}. It also
follows from Theorem \ref{t-2} below, where we give a formula for it. 

\subsection{A local formula for the trace}
Let $\mathcal A_D(M)$ be the quantum algebra with Fedosov connection $D=\nabla+[A,\ ]_\epsilon$.

To give a formula for the trace we use the Lie
algebra cocycle $\Theta_{2n}^N$
obtained from the Hochschild cocycle
$\tau_{2n}$, see \ref{ss-Juno}.
As we consider here only the scalar case, we take $N=1$
in this section and write $\Theta_{2n}$ instead
of $\Theta_{2n}^1$. We will use the following property
of this cocycle:
\begin{enumerate}
\item[I.] $\Theta_{2n}$ extends uniquely to a continuous
$2n$-cocycle of the topological Lie algebra
$\mathcal A_{2n}$ with values in the continuous
dual $(\mathcal A_{2n})^*$.
\item[II.] $\Theta_{2n}$ vanishes if any of its arguments
is in $\mathfrak {sp}_{2n}$ or in the center 
$\Ceps$.
\item[III.]  $\Theta_{2n}$ is $\mathfrak{sp}_{2n}$-invariant: if $a\in \mathfrak{sp}_{2n}$,
\[
\sum_{j=1}^{2n}
\Theta_{2n}(\cdots\wedge[a,a_j]_\epsilon \wedge\cdots)
(f)
+
\Theta_{2n}(a_1\wedge\cdots\wedge a_{2n})
([a,f]_\epsilon )=0,
\]
\item[IV.] $\Theta_{2n}(p_1\wedge q_1\wedge\cdots\wedge
p_n\wedge q_n)(1)=1$.
\end{enumerate}

Properties II and III are the statement that the
cocycle $\Theta_{2n}$ is relative for 
$\mathfrak{gl}_1\oplus\mathfrak{sp}_{2n}$
(Corollary \ref{c-Helena}). Property IV is
\eqref{e-priamus} in Corollary \ref{c-Helena}.

Recall that each fibre of $E(M)$ is isomorphic to
the Weyl algebra and an isomorphism is given for
each choice of symplectic frame of the tangent
bundle. By Property III, $\Theta_{2n}$ is defined 
independently of this choice on 
the fibres of the bundle $E(M)$ and we get a map
$\Theta_{2n}:\wedge^{2n}\Gamma(M,E(M))\otimes \Gamma(M,E(M))\to C^\infty(M)[[\epsilon]]$:
\[
f_1\wedge\cdots\wedge f_{2n}\otimes f\to
\Theta_{2n}(f_1\wedge\cdots\wedge f_{2n})(f).
\]
This map extends by $\Omega^\Cdot(M)$-multilinearity to
a map $\wedge^{2n}\Omega^\Cdot(M,E(M))\otimes 
\Omega^\Cdot(M,E(M))\to \Omega^\Cdot(M)[[\epsilon]]$.

We define a map 
$\psi_D:\mathcal A_D(M)\to \Omega^{2n}(M)[[\epsilon]]$.
\[
\psi_D(f)=\frac{(-1)^n}{(2n)!}\Theta_{2n}(A\wedge\cdots\wedge A)(f).
\]
This map only depends on $D$ and not on the choice of
representation of $D$ as $\nabla+[A,\ ]_\epsilon$. Indeed,
any two representations differ by
changing the symplectic connection, 
which amounts to adding an $\mathfrak{sp}_{2n}$-valued 
1-form to $A$, or adding a $1$-form with values in the
center. By Property II, $\psi_D$ is insensitive to these
changes.

\begin{proposition}\label{p-1}\ \begin{enumerate}
\item[(i)] For any horizontal sections $f,g\in\mathcal A_D(M)$,
$\psi_D([f,g]_\epsilon )\in d\Omega^{2n-1}.$
\item[(ii)] If $D$, $D'$ are gauge equivalent Fedosov connections
and  $\phi:\mathcal A_D(M)\to\mathcal A_{D'}(M)$ 
is the corresponding isomorphism, then for all $f\in 
\mathcal A_D(M)$,
\[
\psi_D(f)-\psi_{D'}(\phi(f))\in d\Omega^{2n-1}.
\]
\end{enumerate}
\end{proposition}

\noindent{\it Proof:} 
The map $\Theta_{2n}:
\wedge^{2n}\Omega^\Cdot(M,E(M))\otimes 
\Omega^\Cdot(M,E(M))\to \Omega^\Cdot(M)[[\epsilon]]$
obeys, by III,
\begin{eqnarray*}
d(\Theta_{2n}(f_1\wedge\dots\wedge f_{2n})(f))
&=&\sum_{i=1}^{2n}\pm\Theta_{2n}(f_1\wedge\cdots\wedge
 \nabla f_i\wedge\cdots\wedge f_{2n})(f)
\\
&&
\pm\Theta_{2n}(f_1\wedge\dots\wedge f_{2n})(\nabla f)
\end{eqnarray*}
(the signs are $(-1)^{\sum_{j=1}^{i-1}\mathrm{deg}(f_j)}$,
$(-1)^{\sum_{j=1}^{2n}\mathrm{deg}(f_j)}$).).

(i)
We claim that $\psi_D([f,g]_\epsilon )=d\varphi$, where
\[
\varphi
=
\frac{(-1)^n}{(2n)!}\Theta_{2n}(A\wedge\cdots\wedge A\wedge f)( g)
\]
The claim follows using the fact that  $\nabla A+\frac12[A,A]_\epsilon =\Omega$
belongs to the center and the identities $\nabla f+[A, f]_\epsilon =
\nabla g+[A, g]_\epsilon =0$:
\begin{eqnarray*}
d\varphi&=&
(2n-1)\Theta_{2n}(\nabla A\wedge A^{2n-2}\wedge f)(g)
\\
&&+\Theta_{2n}(A^{2n-1}\wedge \nabla f)(g)
+\Theta_{2n}(A^{2n-1}\wedge f)(\nabla g)\\
&=&
-(2n-1)\Theta_{2n}(\textstyle{\frac12}
[A,A]_\epsilon\wedge A^{2n-2}\wedge f)(g)
\\
&&-\Theta_{2n}(A^{2n-1}\wedge [A,f]_\epsilon)(g)
-\Theta_{2n}(A^{2n-1}\wedge f)([A,g]_\epsilon)\\
&=&\Theta_{2n}(A^{2n})([f,g]_\epsilon),
\end{eqnarray*}
by the cocycle property of $\Theta_{2n}$.

(ii) A similar calculation show that if $A_t$ is a family of Fedosov
connection forms such that $\frac d{dt}A_t=
\nabla\lambda_t+[A_t,\lambda_t]_\epsilon $
and  $\frac d{dt} f_t=-[\lambda_t, f_t]_\epsilon $, then
$
\frac d{dt}\psi_{D_t}(f_t)=d\varphi_t
$, with 
$\varphi_t=-\frac{(-1)^n}{(2n)!}\Theta_{2n}(A\wedge
\cdots\wedge A\wedge\lambda_t)(f_t).
$
\hfill$\square$
\medskip

\begin{thm}\label{t-2} The unique normalized trace on
$\mathcal A_D(M)$ is given by the formula
\[
\mathrm{Tr}_D(f)=\frac1{(2\pi\ii\epsilon)^n}
\int_M\psi_D(f).
\]
If $D,D'$ are equivalent Fedosov connections and
$\phi$ 
is the corresponding isomorphism, then for all $f\in 
\mathcal A_D(M)$,
\[
\mathrm{Tr}_D(f)=\mathrm{Tr}_{D'}(\phi(f)).
\]
\end{thm}

\noindent{\it Proof:} We know from Prop.~\ref{p-1} (i) that 
the right-hand side is a trace. In  remains to show that it
is normalized. For this we use the fact that locally every 
quantum algebra is isomorphic to the Weyl algebra 
(Lemma \ref{p-0}). By Prop.~\ref{p-1}, (ii), it is sufficient
to evaluate $\psi_D$ for $E(\RR^{2n})$ with the connection
$D^\circ$. In this case, $A=\sum_{i,j}\omega^\circ_{ij}y_idx_j$, so the
claim follows from Property II of $\Theta_{2n}$. The
second claim follows from Prop.~\ref{p-1}, (ii).
\hfill $\square$

\subsection{Remark: non-commutative formal geometry.}
The Fedosov construction has an interpretation in terms of 
a non-commutative version of formal geometry \cite{GK}, which is the
approach that Nest and Tsygan \cite{NT1} take. It turns out that
our form $\psi_D$ also can be naturally interpreted in this context.
For this interpretation it is better to work over the real numbers,
so in this subsection we replace $\mathbb C$ by $\mathbb R$.

We sketch the constructions here, starting with recalling the commutative
case.
Recall that in (commutative)
formal geometry one considers the 
manifold $\tilde M$ of $\infty$-jets at $0$ of local parametrizations
$\mathbb R^d\to M$ of a given $d$-manifold $M$. The linearization at $0$ defines
a $\mathrm{GL}_d$-equivariant projection $\tilde M\to FM$  onto the frame bundle,
with contractible fibers.
The Lie algebra $W_d$ of 
formal vector fields acts freely and transitively on $\tilde M$
by infinitesimal reparametrizations. This means that we have for each
$\varphi\in\tilde M$ 
an isomorphism $W_d\to T_\varphi\tilde M$, whose inverse defines a 
$W_d$-valued
one-form $\omega_\mathrm{MC}$. The action property is equivalent
to the fact that $\omega_\mathrm{MC}$ is an
$H$-equivariant Maurer--Cartan form, namely
that it obeys  the Maurer--Cartan equation 
$d\omega_\mathrm{MC}+\frac12[\omega_\mathrm{MC},\omega_\mathrm{MC}]=0$. Moreover
the action of the Lie subalgebra of linear vector fields can be integrated to
an action of $\mathrm{GL}_d$ by linear coordinate transformation.
 
Abstractly, we thus have a manifold $X=\tilde M$ with a Maurer--Cartan one-form
$\omega_\mathrm{MC}$ and a subalgebra $\mathfrak h$ whose action integrates
to an action of a Lie group $H$ for which $X$ is a principal $H$-bundle.
This gives rise to a canonical ring homomorphism
$H^\Cdot(\mathfrak g,H)\to H^\Cdot(X)_{H-\text{basic}}\simeq
H^\Cdot(X/H)$. This map is 
given by evaluating a cocycle on the Maurer--Cartan form. More generally,
if $V$ is a $\mathfrak g$-module with integrable action of $H$, then
the associated bundle $E_V=X\times_{H}V\to X/H$ carries a flat connection $D$ 
coming from the connection $d+\omega_\mathrm{MC}$ on the trivial
bundle $X\times V$. We then obtain a canonical homomorphism
$H^\Cdot(\mathfrak g,H;V)\to H^\Cdot(X/H,E_V)$ from the cohomology
with coefficients in $V$ to the cohomology of the complex of $E_V$-valued
differential forms with differential $D$. In view of the non-commutative
generalization, suppose we have a central extension
$0\to \mathfrak z\to\hat{\mathfrak g}\to\mathfrak g$ and a lift 
$\hat\omega_\mathrm{MC}$ of the Maurer--Cartan form to a $\hat g$-valued
1-form. Then we get the Maurer--Cartan equation in the form
\begin{equation}\label{e-MChat}
d\hat\omega_\mathrm{MC}+\frac12[\hat\omega_\mathrm{MC},\hat\omega_\mathrm{MC}]=F,
\end{equation}
for some $H$-basic closed 2-form $F$ with values in the center, whose cohomology class
in $H^\Cdot(X/H,\mathfrak z)$ is independent of the choice of the lift.

Turning to the non-commutative case,
suppose that $M$ is a symplectic manifold with a star-product. 
It follows from the local nature of star-products that the
space of $\infty$-jets $\mathrm{Jet}_x(M)$ of
functions in $C^\infty(M)[[\epsilon]]$
at any point $x\in M$
also comes with a star-product. Moreover locally
any star-product on a symplectic manifold is
equivalent to the Moyal product. 
We then define $\tilde M$ to be
the space of pairs $(x,\varphi^*)$, where $x\in M$ and $\varphi^*$
is an isomorphism $\mathrm{Jet}_x(M)\to\mathcal A^+_{2n}$.
Passing to 1-jets, we get an $\Sp_{2n}$-equivariant
projection $p:\tilde M\to F_{\Sp}M$ onto the symplectic
frame bundle.
This projection has contractible fibers: any two isomorphisms differ by an 
automorphism of the Weyl algebra. Such an
automorphism is given on topological generators by
$y_j\mapsto f_j(y,\epsilon)=\sum_k A_{jk}y_k+\cdots$ obeying $[f_j,f_k]_\epsilon=
\delta_{jk}$, for some symplectic matrix $A$. A contracting
homotopy
is given by $t^{-1}f_j(ty,t^2\epsilon)$, $t\in[0,1]$.
The Lie algebra $\mathrm{Der}\, \mathcal A^+_{2n}$ acts on
$\tilde M$ and we get a Maurer--Cartan 1-form $\omega_\mathrm{MC}$ with values in $\mathrm{Der}\, \mathcal A^+_{2n}$. 
Lifting $\omega_\mathrm{MC}$ to a form
$\hat\omega_\mathrm{MC}$ with values in $\mathcal A^+_{2n}$, we get the Maurer--Cartan equation in the form
\eqref{e-MChat}
for some $\Sp_{2n}$-basic closed 2-form $F$ with values in $\mathfrak z=\mathbb R[[\epsilon]]$. 
The cohomology class of $F$ is independent of 
the choice of lift. Its pull-back by an $\Sp_{2n}$-equivariant section of
$\tilde M\to F_\Sp M$ gives a basic class in $H^2(F_\Sp M,
\mathbb R[[\epsilon]])$
and thus a class in $H^2(M,\mathbb R[[\epsilon]])$.
This is the characteristic class of the star-product. 

A Fedosov quantum algebra $\mathcal A^+_D(M)$ comes with a section $\sigma:F_{\Sp}M\to \tilde M$ of $p$.
It sends a symplectic frame $(x,e_1,\dots,e_{2n})$ 
at $x\in M$ to
$(x,\varphi)$ with $\varphi$ constructed as follows: let
$[f]\in \mathrm{Jet}_x(M)$ be represented by a function
$f$ defined in some neighborhood $U$ 
of $x$ and let $\hat f$
be the corresponding horizontal section of 
$\mathcal A_D^+(U)$. Then $\varphi([f])$ is the
value at $x$ of $\hat f$ via the identification of the
fibre at $x$ with $\mathcal A^+_{2n}$ given by
the symplectic frame. It is easy to see that 
this gives a well-defined map on jets.

Again we get a homomorphism $H^\Cdot(\mathrm{Der}\,
\mathcal A^+_{2n},\mathfrak{sp}_{2n}; V)\to H^\Cdot(M)$,
for any $\mathrm{Der}\,\mathcal A^+_{2n}$-module $V$ (since $\Sp_{2n}$ is connected we may
pass to the Lie algebra). Equivalently, considering $V$ as a module over the Lie algebra
$\mathcal A^+_{2n}$ with trivial action of the center, we have a formal geometry homomorphism
\[
 H^\Cdot(\mathcal A^+_{2n},\mathfrak{sp}_{2n}\oplus \mathbb R[[\epsilon]]; V)
\to H^\Cdot(\tilde M/\Sp_{2n},E_V).
\]
If $V=\mathcal A^+_{2n}$, $H^0(\tilde M/\Sp_{2n},E_V)=\mathcal A^+_D(M)$ is the Fedosov quantum
algebra (the pull-back by $\sigma$ of $E_V$ is the Weyl algebra bundle $E(M)$). 
If $V=(\mathcal A^+_{2n})^*$ is the continuous dual $\mathbb R[[\epsilon]]$-module,
with coadjoint action of $\mathcal A^+_{2n}$,
we obtain in particular a map
\begin{eqnarray*}
\lefteqn{ H^{2n}(\mathcal A^+_{2n},\mathfrak{sp}_{2n}\oplus \mathbb R[[\epsilon]]; (\mathcal A^+_{2n})^*)
\to 
H^{2n}(\tilde M/\Sp_{2n},E_{(\mathcal A^+_{2n})^*})}\\
&
\to& 
\mathrm{Hom}_{\mathbb R[[\epsilon]]}\left(H^0(\tilde M/\Sp_{2n}
,E_{\mathcal A^+_{2n}}),H^{2n}(M,\mathbb R[[\epsilon]])\right).
\end{eqnarray*}
The second arrow is induced by the pairing between $\mathcal A^+_{2n}$ and its dual.
Up to the normalization factor $(2\pi\mathrm i\,\epsilon)^{-n}$ and passing to 
Laurent series with complex coefficients, 
the image of our basic class $\Theta_{2n}$ by this
composition is the map $\psi_D:\mathcal A_D(M)\to H^{2n}(M,\Ceps)$ 
on the level of cohomology.
\section{The local Riemann--Roch--Hirzebruch theorem}
\label{sect4}
In this section, we study the evaluation of the relative
$2n$-cocycle 
$\Theta^N_{2n}\in 
C^{2n}(\mathfrak{gl}_N(\mathcal A^\pol_{2n}),\mathfrak{sp}_{2n}
\oplus \mathfrak{gl}_N;\mathfrak{gl}_N(\mathcal A^\pol_{2n})^*)$
introduced in Section \ref{sect2}
at the identity operator and relate it to characteristic
classes. We obtain a direct proof of
the local algebraic Riemann--Roch formula
discovered in \cite{FT2}

\subsection{Chern--Weil homomorphism}\label{s-51}
Let $\mathfrak g$ be a Lie algebra and $\mathfrak h\subset
\mathfrak g$ be a subalgebra. Suppose that there is
an $\mathfrak h$-equivariant projection $\pr:\mathfrak g
\to \mathfrak h$, that is a map commuting with the
adjoint action of $\mathfrak h$ and restricting to
the identity on $\mathfrak h$. 
The amount by which $\pr$ fails to be a Lie 
algebra homomorphism is measured by the 
``curvature'' $\Curv\in\mathrm{Hom}(\wedge^2\mathfrak g,
\mathfrak h)$:
\[
\Curv(v\wedge w)=[\pr(v),\pr(w)]-\pr([v,w]).
\]
Then there is 
a Chern--Weil graded algebra homomorphism
\[
\chi: (S^\Cdot\mathfrak h)^{*\mathfrak h}
\to H^{{}^2\Cdot}(\mathfrak g,\mathfrak h),
\]
independent of the choice of projection $\pr$,
sending $P\in (S^q\mathfrak h)^{*\mathfrak h}$ 
to the class of $\chi(P)=P(C^q)$. Explicitly,
\begin{eqnarray*}
\lefteqn{\chi(P)(v_1\wedge\cdots\wedge v_{2q})}\\
&&=\frac1{q!}\sum_{
\sigma
}{}'\,\mathrm{sign}(\sigma)
P\left(
\Curv(v_{\sigma(1)},v_{\sigma(2)}),
\dots,
\Curv(v_{\sigma(2q-1)},v_{\sigma(2q)})\right).
\end{eqnarray*}
The sum is over all permutations $\sigma\in S_{2q}$
such that $\sigma(2i-1)<\sigma(2i)$. We give a
review of this construction below.

In our case, we take $\mathfrak g=\mathfrak{gl}_N(\mathfrak A^\pol_{2n})$, $\mathfrak h=\mathfrak{gl}_N\oplus\mathfrak{sp}_{2n}$ with $\mathfrak{sp}_{2n}$ embedded as $1\otimes 
\mathfrak{sp}_{2n}$
and $\mathfrak{gl}_N$ as $\mathfrak{gl}_N\otimes 1$.
Then we have a projection $\pr=\pr_1\oplus \pr_2$ with
$\pr_1(M\otimes a)=\frac1N\mathrm{tr}(M)\,a_2\in\mathfrak{sp}_{2n}$ and $\pr_2(M\otimes a)=Ma_0$, where
$a_j$ is the component of $a\in K[y_1,\dots,y_n]$ 
homogeneous of degree
$j$ in $y$.

\subsection{Characteristic classes}

We introduce ad-invariant functions on $\mathfrak{sp}_{2n}$
and $\mathfrak{gl}_N$
(the A-roof genus and the Chern character). An invariant 
 function 
on $\mathfrak{sp}_{2n}$ is uniquely determined by its
values on diagonal symplectic
matrices $\mathrm{diag}(t_1,-t_1,\dots,t_n,-t_n)$ and
a polynomial function of $t_1,\dots,t_n$ extends to
an invariant function if and only if it is invariant
under the Weyl group, which is generated by
permutations and sign reversals $t_i\to -t_i$. To
the function $\prod_{i=1}^n\sinh(t_i/2)^{-1}t_i/2$ there
corresponds the A-roof genus
\begin{eqnarray*}
\hat{\mathrm A}(X)&=&\mathrm{det}\left(\frac {X/2} {
\sinh(X/2)}\right)^{\frac12}\\
&=&1-\frac1{48}\mathrm{tr}(X^2)
+\frac1{4608}\mathrm{tr}(X^4)
+\frac1{5760}(\mathrm{tr}(X^2))^2+\cdots.
\end{eqnarray*}
We will
need the rescaled version $\hat {\mathrm A}_\epsilon(X)=
\hat{\mathrm A}(\epsilon X)$.
The Chern character
\[
\mathrm{Ch}(Y)=\mathrm{tr}\left(\exp Y\right),
\]
is an ad-invariant function of $Y\in \mathfrak{gl}_N$.
These functions are analytic at $0$. The Taylor expansion at zero 
of their product is a sum of invariant polynomials.
If $X=X_1\oplus X_2\in\mathfrak{sp}_{2n}\oplus \mathfrak{gl}_N$, 
denote $(\hat{\mathrm{A}}_\epsilon\mathrm{Ch})_j\in 
(S^{j}\mathfrak h^*)^{\mathfrak h}$
the  term homogeneous of degree $j$
of the Taylor expansion of $(\hat{\mathrm A}_\epsilon\,\mathrm{Ch})(X)=
\hat{\mathrm A}_\epsilon(X_1)\mathrm{Ch}(X_2)$: it is
defined a symmetric $j$-linear form by
\[
\hat{\mathrm{A}}_\epsilon\mathrm{Ch}(X)
=1+(\hat{\mathrm{A}}_\epsilon\mathrm{Ch})_1(X)
+\frac1{2!}(\hat{\mathrm{A}}_\epsilon\mathrm{Ch})_2(X,X)
+\cdots
\]

\subsection{Statement of the Theorem}\label{ss-MrsLucreziaWarrenSmith}

Let the morphism of $\mathfrak{gl}_N(\mathcal A^\pol_{2n})$-modules
$\ev:\mathfrak{gl}_N(\mathcal A^\pol_{2n})^*\to \Ceps$ 
be the evaluation at the identity. Then $\ev\Theta_{2n}$ is a cocycle
in $C^{2n}(\mathfrak{gl}_N(\mathcal A^\pol_{2n}),\mathfrak{sp}_{2n}\oplus \mathfrak{gl}_N;\Ceps)$.
It is given by the formula
\[
\ev\Theta_{2n}(v_1,\dots,v_{2n})=\Theta_{2n}(v_1,\dots,v_{2n})(1).
\]
Let $\mathfrak g=\mathfrak{gl}_N(\mathcal A^\pol_{2n})$,
$\mathfrak h=\mathfrak{sp}_{2n}\oplus\mathfrak{gl}_N$.

\begin{thm}\label{t-lRRH}
$
[\ev\Theta_{2n}]=
{(-1)^n}\chi((\hat{\mathrm{A}}_\epsilon\,\mathrm{Ch})_n).
$
\end{thm}

The proof of this theorem consists of evaluating explicitly
the integral defining $\Theta_{2n}$ on special arguments. Before
giving this proof we need to introduce some technology to 
show that it is indeed sufficient to prove the formula on
these special arguments. 

\subsection{The relative Weil algebra} Let
$\mathfrak g$ be a Lie algebra, 
$\mathfrak h\subset \mathfrak g$
a finite dimensional subalgebra, $\pi:\mathfrak g\to
\mathfrak g/\mathfrak h$ the canonical projection. 
 The relative
Weil algebra is the differential graded
commutative algebra $W(\mathfrak g,\mathfrak h)=\oplus_{i,j\geq0} W^{i,j}(\mathfrak g,\mathfrak h)$,
where $W^{i,j}(\mathfrak g,\mathfrak h)=(\wedge^i(\mathfrak g/\mathfrak h)
\otimes S^j\mathfrak g)^{*\mathfrak h}$ has grading
$i+2j$. The differential $d=d_{\mathrm{Lie}}+d'$
has two pieces: $d_{\mathrm{Lie}}:W^{i,j}\to W^{i+1,j}$ is
the differential of the relative Lie algebra cochain
complex with coefficients in the $\mathfrak g$-module $S^{j}\mathfrak g^*$
and $d':W^{i,j}\to W^{i-1,j+1}$ is the map
\[
d'c(v\otimes w_1\cdots w_{j+1})=\sum_{k=1}^{j+1}
c(\pi(w_k)\wedge v\otimes
w_1\cdots\hat w_k\cdots w_{j+1}), 
\]
$v\in\wedge^{i-1}\mathfrak g/\mathfrak h, w_k\in\mathfrak g$. The product in $W(\mathfrak g,\mathfrak h)$ is
given by a sum over shuffles (with proper signs)
as in the case of the Lie algebra cochain complex,
see \ref{ss-HLAC}.
We have $C^\Cdot(\mathfrak g,\mathfrak h)=W^{\Cdot,0}(
\mathfrak g,\mathfrak h)$ with differential $d_\mathrm{Lie}$. This differential graded commutative algebra can
thus be identified with $W(\mathfrak g,\mathfrak h)
/\oplus_{j\geq 1}W^{\Cdot,j}(\mathfrak g,\mathfrak h)$.
So we have a canonical projection
of differential graded commutative algebras 
$W^\Cdot(\mathfrak g,\mathfrak h)\to 
C^\Cdot(\mathfrak g,\mathfrak h)$.

If $\mathfrak h=0$, the Weil algebra is acyclic \cite{C}.
More generally one has a projection $W(\mathfrak g,
\mathfrak h)\to W(\mathfrak h,\mathfrak h)=(S^\Cdot\mathfrak h)^{*\mathfrak h}$ which induces a homomorphism
on cohomology, and one has the following well-known result.

\begin{lemma} \label{l-WA} Suppose there exists
an $\mathrm{ad}\,\mathfrak h$-invariant  subspace $V\subset \mathfrak g$ so that
$\mathfrak g=\mathfrak h\oplus V$.
Then $H^{2j+1}(W(\mathfrak g,\mathfrak h))=0$ and 
$H^{2j}(W(\mathfrak g,\mathfrak h))\simeq(S^{j}
\mathfrak h)^{*\mathfrak h}$,
$(j=0,1,2,\dots)$.
\end{lemma}

\medskip

\noindent{\it Proof:}
Consider the filtration $W(\mathfrak g,\mathfrak h)
=F^0\supset F^1\supset\cdots$, with $F^p=
\oplus_{i+j\geq p}W^{i,j}(\mathfrak g,\mathfrak h)$. In
the resulting spectral sequence $E_0^{p,q}=
W^{p-q,q}(\mathfrak g,\mathfrak h)$ with differential $d'$.
Let $V$ be an $\mathrm{ad}\,\mathfrak h$-invariant complement.
Then $E_0\simeq(\wedge^\Cdot V\otimes S^\Cdot V
\otimes S^\Cdot\mathfrak h)^{*\mathfrak h}$.
Introduce a derivation $h$ 
of degree $-1$ on $E_0$ by 
\[
hc(v_1\wedge\cdots\wedge v_i\otimes w\otimes u)
=\sum_{k=1}^i(-1)^{k-1}
c(v_1\wedge\cdots\hat v_k\dots\wedge v_i\otimes v_kw
\otimes u),
\]
Then for
$a\in (\wedge^i V\otimes 
S^jV\otimes S^l\mathfrak h)^{*\mathfrak h}$,
\[
(d'h+hd')(a)=(i+j)a,
\]
and every cocycle with $i+j>0$ is a coboundary.
It follows that 
$E_1^{p,p}=(S^p\mathfrak h)^{*\mathfrak h}$,
$E_1^{p,q}=0$, $q\neq p$. The spectral sequence degenerates
and the result follows. \hfill $\square$

\medskip

We now give a formula for a cocycle 
corresponding to a given element
of $(S^j\mathfrak h)^{*\mathfrak h}$ using a
version of the Chern--Weil construction. The
algebra $W(\mathfrak g,\mathfrak h)$ consists
of {\em basic} elements in the absolute Weil algebra
$W(\mathfrak g)=W(\mathfrak g,0)$. Namely, for
each $h\in\mathfrak h$, we have
the internal multiplication $\iota(h):
W^{i,j}(\mathfrak g)
\to
W^{i-1,j}(\mathfrak g)$, 
such that $\iota(h)c(v\otimes w)=c(h\wedge v\otimes w)$
and the action of $h\in\mathfrak h$ on $W(\mathfrak g)$
is given by the Cartan formula $L(h)=d\circ \iota(h)+
\iota(h)\circ d$. Basic elements are those in the
kernel of $i(h), L(h)$ for all $h\in\mathfrak h$.
Let
$V$ be an $\mathrm{ad}\,\mathfrak h$-invariant complement
of $\mathfrak h$ in $\mathfrak g$.
Then the projection $\pr:\mathfrak g
\to \mathfrak h$ along $V$ is $\mathfrak h$-equivariant.
We may regard it as an element $A$ 
of $\wedge^1\mathfrak g^*
\otimes
\mathfrak h\subset W(\mathfrak g)\otimes \mathfrak h$.
The equivariance and projection property of $\pr$
imply that $A$ is a {\em connection} in the sense of
\cite{C}. Namely, an element of degree 1 obeying
\[
(\iota(h)\otimes \mathrm{id})(A)=h,\qquad
(L(h)\otimes \mathrm{id}+\mathrm{id}\otimes \mathrm{ad}\,h)
(A)=0.
\]
Because of this, we have a characteristic
homomorphism of  differential graded commutative
algebras\[
\chi_W:W(\mathfrak h,\mathfrak h)=
(S\mathfrak h)^{*\mathfrak h}
\to W(\mathfrak g,\mathfrak h).
\]
It is constructed as follows. The curvature
$F=dA+\frac12[A,A]$
of $A$ in the differential graded Lie algebra 
$W(\mathfrak g)\otimes \mathfrak h$
is $\mathfrak h$-equivariant and obeys
$(\iota(h)\otimes\mathrm{id})(F)=0$. It follows that
$F^j\in W(\mathfrak g)\otimes S^j\mathfrak h$ has
the same properties. Thus
\[
\chi_W(P)=\frac1{j!}(\mathrm{id}\otimes P)(F^j),\qquad 
P\in (S^j\mathfrak h)^{*\mathfrak h}
\]
is basic. Passing to cohomology, we have a Chern--Weil
homomorphism
\[
\chi_W: (S \mathfrak h)^{*\mathfrak h}
\to H(W(\mathfrak g,
\mathfrak h)),
\]
which is independent of the choice of
$\mathfrak h$-equivariant projection.

\begin{proposition} \label{p-Venus}
Assume that $\mathfrak h$ has an 
$\mathrm{ad}\,\mathfrak h$-invariant 
complement in $\mathfrak g$.
Then the map $\chi_W$ is an isomorphism.
Its composition $\chi:(S\mathfrak h)^{*\mathfrak h}\to
H(W(\mathfrak g,\mathfrak h))\to H(\mathfrak g,
\mathfrak h)$ with the canonical map (induced by
the projection $W(\mathfrak g,\mathfrak h)\to C(\mathfrak g,\mathfrak h)$) is the Chern--Weil homorphism: it sends
$P\in (S^{j}\mathfrak h)^{*\mathfrak h}$ 
to the class of 
\begin{equation}\label{e-Aeolus}
\chi(P)=P(F_1^j)/j!,\qquad 
F_1(v,w)=[\pr(v),\pr(w)]
-\pr([v,w]),
\end{equation}
$v,w\in\mathfrak g$,
  for any choice of 
 $\mathrm{ad}\,\mathfrak h$-equivariant projection
$\pr:\mathfrak g
\to \mathfrak h$.
\end{proposition}

\noindent{\it Proof:}
The curvature  $F=dA+\frac12[A,A]$ 
has two components: $F=F_1+F_2$, $F_1=d_{\mathrm{Lie}}A
+\frac12[A,A]\in W^{2,0}(\mathfrak g)\otimes\mathfrak h$,
$F_2=d'A\in W^{0,1}(\mathfrak g)\otimes \mathfrak h$.
Identifying $W(\mathfrak g)\otimes\mathfrak h$ with
$\mathrm{Hom}(\wedge \mathfrak g\otimes S\mathfrak g,\mathfrak h)$, we have
 $F_1(v,w)=d_{\mathrm{Lie}}A(v,w)+[A(v),A(w)]$
and the formula \eqref{e-Aeolus}
for $F_1$ follows. For $v\in S^1\mathfrak g$,
$F_2(v)=\pr(v)$. These formulae allow us to
compute the components of $\chi_W(P)$. If 
$P\in (S^j\mathfrak h)^{*\mathfrak h}$, $\chi_W(P)=\sum_{k=0}^j 
\chi^k_W(P)$, with $\chi_W^k(P)\in W^{2k,j-k}(\mathfrak g,
\mathfrak h)$. The
component in $W^{0,j}=(S^j\mathfrak g)^{*\mathfrak h}$ is
\[
\chi^0_W(P)(v_1\cdots v_j)=
P(\pr(v_1)\cdots \pr(v_j)).
\]
This shows that $\chi_W$ is a right inverse  to
the restriction map $W(\mathfrak g,\mathfrak h)\to
W(\mathfrak h,\mathfrak h)=(S\mathfrak h)^{*\mathfrak h}$,
which by Lemma \ref{l-WA} is an isomorphism. 
Thus $\chi_W$ is an isomorphism.

The component of $\chi_W(P)$ in $W^{2j,0}$ is
$\chi_W^j(P)=P(F_1^j)/j!$. This proves the second
part of the claim. \hfill $\square$

\subsection{Characteristic homomorphism for the Weyl
algebra}

\begin{proposition}\label{p-Achates}
 Let $\mathfrak g=\mathfrak{gl}_N
(\A^\pol_{2n})$, $\mathfrak h=\mathfrak{gl}_N\oplus 
\mathfrak{sp}_{2n}$, $N>>n$. Then for all $j=0,\dots,n$,
the Chern--Weil homomorphism 
$\chi:(S^{j}\mathfrak h^*)^\mathfrak h\to H^{2j}
(\mathfrak g,\mathfrak h)$
is an isomorphism.
\end{proposition}

\noindent 
To prove this proposition we need a statement analogous to
Theorem \ref{t-Dido} for relative cohomology.

\begin{lemma}\label{l-Ilionius}
Let $j\geq 0$. If $k\leq 2n$,
$H^k(\mathfrak g,\mathfrak h; S^j\mathfrak
g^*)=H^k(\mathfrak g;S^j\mathfrak g^*)$. In
particular,
$H^k(\mathfrak g,\mathfrak h; S^j\mathfrak g^*)=0$
if $k<2n$ and $j\geq1$.
\end{lemma}

\noindent{\it Proof:}
The Hochschild--Serre spectral sequence
\cite{HS}
for the pair $(\mathfrak g ,\mathfrak h)$ converging
to $H^\Cdot(\mathfrak g;S^j\mathfrak g^*)$ has $E_1$-term
$E_1^{p,q}=H^q(\mathfrak h;
\mathrm{Hom}(\wedge^p\mathfrak g/\mathfrak h,
S^j\mathfrak g^*))$. 
The 
$\mathfrak h$-module 
$\mathrm{Hom}(\wedge^p\mathfrak g/\mathfrak h,
S^j\mathfrak g^*)$ is
semisimple. For any finite-dimen\-sional Lie algebra $\mathfrak h$ 
with a non-degenerate invariant symmetric 
bilinear form and
semisimple $\mathfrak h$-module $M$, we have $H^q(\mathfrak h;M)=
H^q(\mathfrak h)\otimes M^\mathfrak h$.
[Proof: the operator on cochains given by
composition with the
quadratic Casimir operator $Q=\sum e_ie^i$ associated
to the bilinear form is homotopic
to zero: 
$Q\circ c=d h(c)+h (dc)$, 
with homotopy $h(c)(a_1,\dots,a_{p-1})=\sum e^i
c(e_i,a_1,\dots,a_{p-1})$, $c\in C^p(\mathfrak h,M)$].
In our case, we obtain
$E_1^{p,q}=H^q(\mathfrak h)\otimes
\mathrm{Hom}
(\wedge^p(\mathfrak g/\mathfrak h)^*,S^j\mathfrak g^*)^\mathfrak h$. The differential
is the relative cochain differential on the second
factor. Hence
\[
E_2^{p,q}=H^q(\mathfrak h)\otimes H^p(\mathfrak g,\mathfrak h;S^j\mathfrak g^*).
\]
Now we use the fact that the spectral sequence converges
to $H^{p+q}(\mathfrak g,S^j\mathfrak g^*)=0$ 
if $p+q<2n$ and $j>0$:
we have $E^{0,0}_2=0$ and, since $H^0(\mathfrak h)=\CC$,
$H^0(\mathfrak g,\mathfrak h;S^j\mathfrak g^*)=0$.
If we assume inductively that 
$H^p(\mathfrak g,\mathfrak h;S^j\mathfrak g^*)=0$ 
for $p=0,\dots, p_0<2n$, we
have $E_2^{p,q}=0$ for all $p<p_0$ and all $q$. It
follows that $E_2^{p_0,0}=H^{p_0}(\mathfrak g;S^j\mathfrak g^*)=0$ and thus $H^{p_0}(\mathfrak g,\mathfrak h;S^j\mathfrak g^*)=0$, proving the induction step and the claim
for $k<2n$. Since $E_2^{p,q}=0$ for all $p<2n$ and all 
$q$, the map $E_{2}^{2n,0}=H^{2n}(\mathfrak g,\mathfrak h;
S^j\mathfrak g^*)\to H^{2n}(\mathfrak g,S^j\mathfrak g^*)$
is an isomorphism.
\hfill $\square$

\medskip

\noindent{\it Proof of Proposition \ref{p-Achates}:}
Let $E_r^{p,q}$ be the spectral sequence
corresponding to the filtration 
$F^p(W(\mathfrak g,\mathfrak h))=\oplus_{j\geq p}
 W^{i,j}(\mathfrak g,\mathfrak h)$ of the Weil algebra (it is not the
same filtration as in the proof of Lemma \ref{l-WA}). 
Then
$E_0^{p,q}=C^q(\mathfrak g,\mathfrak h;S^p\mathfrak g^*)$ is the relative Lie
algebra cochain complex (the second piece $d'$ of
the Weil differential does not contribute). Thus
$E_1^{p,q}=H^q(\mathfrak g,\mathfrak h;S^p\mathfrak g^*)
$ and the canonical 
homomorphism $H(W(\mathfrak g,\mathfrak h))\to
H(\mathfrak g,\mathfrak h)$ is the edge homomorphism
to $E_1^{0,q}$.
By Lemma \ref{l-Ilionius} $E_1^{p,q}=0$  for $p\geq1$ and
$q<2n$. Therefore the edge homomorphism 
is an isomorphism for all $q\leq 2n$. This implies
that  $H^q(\mathfrak g,\mathfrak h)=0$ for $q<2n$ odd
and that
$\chi:(S^{\frac q2} \mathfrak h^*)^\mathfrak h\to H^{q}(\mathfrak g,\mathfrak h)$ is an isomorphism for $q\leq 2n$ even.
\hfill $\square$

\subsection{Proof of Theorem \ref{t-lRRH}}

Let $W_{n,N}$ be the Lie subalgebra of $\mathfrak{gl}_N(\A^\pol_{2n})$
consisting of elements of the form
$
\sum_i f_i(q)p_i\otimes 1+ \sum_j g_j(q)\otimes M_j
$
with
$f_i(q),g_j(q)\in K[q_1,\dots,q_n]$ and $ M_j\in M_N(\CC)$.
It is 
isomorphic to the semidirect product
of the Lie algebra of polynomial
vector fields (derivations
of $K[q_1,\dots,q_n]$) by $\mathfrak{gl}_N(K[q_1,\dots,q_n])$.
Set as before $\mathfrak h=\mathfrak{gl}_N\oplus
\mathfrak{sp}_{2n}$. The Lie algebra $\mathfrak h_1=
W_{n,N}\cap \mathfrak h$ is isomorphic to $\mathfrak{gl}_N
\oplus \mathfrak{gl}_{n}$.
Consider the commutative square
\begin{equation}\label{e-MrsDalloway}
\begin{array}{ccc}
(S^{n}\mathfrak h)^{*\mathfrak h}&\longrightarrow&
(S^{n}\mathfrak h_1)^{*\mathfrak h_1}\\
\downarrow& &\downarrow\\
H^{2n}(\mathfrak{gl}_N(\A^\pol_{2n}),\mathfrak h)
&\longrightarrow&
H^{2n}(W_{n,N},\mathfrak h_1)
\end{array}
\end{equation}
The vertical arrows are Chern--Weil homomorphisms and the
horizontal arrows are induced by the restriction maps.
The relative cohomology of $W_n$ has been calculated by
Gelfand and Fuchs, see \cite{Fu}, using invariant
theory; the same method
can be applied to the case of $W_{n,N}$. The result
can be formulated by saying 
that $C^{k}(W_{n,N},\mathfrak h_1)=0$ for odd $k$
and the Chern--Weil homomorphism
$(S^j\mathfrak h_1)^{*\mathfrak h_1}
\to C^{2j}(W_{n,N},\mathfrak h_1)$ is an isomorphism
for $j\leq n$ and is zero for $j>n$.
In particular, the differential vanishes on 
$C^{\Cdot}(W_{n,N},\mathfrak h_1)$. Thus the vertical
arrows in \eqref{e-MrsDalloway}
are isomorphisms. The upper horizontal arrow (and so
also the lower one) is injective, since invariant polynomials are uniquely determined by their value on the
Cartan subalgebra spanned by $q_ip_i\otimes 1,
1\otimes D$, with $D$ diagonal, for both $\mathfrak h$
and $\mathfrak h_1$. Therefore to determine the
polynomial $P$ such that $\iota_{n}(P)=\ev\Theta_{2n}$
it is sufficient to consider the restriction of $\ev\Theta_{2n}$ to $W_{n,N}$.

We next calculate explicitly the polynomial $P$ so that 
$[\ev\Theta_{2n}]=\chi(P)$. As $P$ is ad-invariant, it
is determined by its value on the Cartan subalgebra
 $\mathfrak d$
spanned by $1\otimes p_iq_i$, $1\leq i\leq n$, and 
$E_{rr}\otimes 1$ ($E_{rr}$, $1\leq r\leq N$
 is the standard basis of the space of diagonal matrices).

\dontprint{ As $P_n$ is $\mathfrak h_1$-invariant, it is sufficient
to consider the restriction of  $P_n$ to diagonal matrices.
To get the signs right, notice that
the Lie algebra $\mathfrak{gl}_n$ is embedded in $W_n$ via $(a_{ij})\mapsto
-\sum a_{ij}x_i\partial/\partial x_j$ and $\mathfrak{gl}_N$ is embedded
as the subalgebra of constant matrices in $\mathfrak{gl}_N(\mathcal O_n)$. Thus a
diagonal matrix $(\mathrm{diag}(t_1,\dots,t_n),\mathrm{diag}(s_1,\dots,s_N))$
corresponds to $-\sum t_ix_i\partial_i+\sum s_r E_{rr}\in 
W_n\ltimes \mathfrak{gl}_N(\mathcal O_n)$. 
}

Introduce $P_n\in (S^n\mathfrak h)^{*\mathfrak h}$
so that its restriction to $\mathfrak d$ is
\begin{eqnarray}
\lefteqn{P_n(M_1\otimes a_1,\dots,M_n\otimes a_n)
=\mathrm{tr}(M_1\cdots M_n)}\notag\\
&&
\times \mu_n\ \int_{[0,1]^n}
\prod_{1\leq i\le j\leq n} \label{e-Purvis}
e^{\epsilon\psi(u_i-u_j)\alpha_{ij}}
(a_1\otimes\cdots\otimes a_n)du_1\cdots du_n,
\end{eqnarray}
with $\mu_n(a_1\otimes\cdots\otimes a_n)=\prod_{i=1}^na_i(0)$.

\begin{lemma} \label{l-525} 
$[\ev\Theta_{2n}]={(-1)^n}\chi(P_n)$. 
\end{lemma}

\noindent{
\it Proof:} As the horizontal restriction maps in
\eqref{e-MrsDalloway} are injective, it is sufficient
to prove the restriction of this identity to $W_{n,N}$.
Since $H^{2n}(W_{n,N},\mathfrak h_1)=C^{2n}(W_{n,N},\mathfrak
h_1)$ the identity becomes an identity of cocycles,
rather than  of cohomology classes.

As the right vertical arrow in \eqref{e-MrsDalloway} is
an isomorphism, we know that there exists a $Q_n\in
(S^n\mathfrak h_1)^{*\mathfrak h_1}$ such that
$\chi(Q_n)$ is the restriction to $W_{N,n}$ of
$[\ev\Theta_{2n}]$. This means that
for all $a_1,\dots, a_{2n}\in W_{n,N}$,
\begin{eqnarray}\label{e-SeptimusWarrenSmith}
\lefteqn{
\ev\Theta_{2n}(a_1\wedge\cdots\wedge a_{2n})
}\notag \\
&=&
\frac1{n!}\sum_{\sigma}{}'\,\mathrm{sign}(\sigma)\,
Q_n(C(a_{\sigma(1)},a_{\sigma(2)}),
\dots,C(a_{\sigma(2n-1)},a_{\sigma(2n)})),
\end{eqnarray}
where the sum is over permutations $\sigma\in S_{2n}$
so that $\sigma(2j-1)<\sigma(2j)$.

We have to show that $Q_n=(-1)^nP_n$. 
By $\mathfrak h$-invariance
it is sufficient to show this on
the Cartan subalgebra $\mathfrak d$
We can obtain the value of $Q_{n}$ on $\mathfrak d$
by evaluating the cocycle $\ev \Theta_{2n}$ on special
vectors. We construct these special vectors out of the
elements $u_{ij},v_{ir}\in W_{n,N}$
$(i,j=1\dots n, \ r=1\dots N)$:
\[
u_{ij}=\left\{\begin{array}{ll}
 \frac12{q_i^2p_i}, & i=j,\\
q_iq_jp_j,& i\neq j,
\end{array}\right. \qquad
v_{ir}=E_{rr}q_i.
\] 
(We write simply $Mf(p,q)$ instead of $M\otimes f(p,q)$).
 These vectors are in the kernel of $\pr$ and obey
\begin{equation}\label{e-PeterWalsh}
[p_i,u_{ij}]_\epsilon=\frac{\partial u_{ij}}{\partial q_i}=
q_jp_j,\qquad [p_i,v_{ir}]_\epsilon
=\frac{\partial v_{ir}}{\partial q_i}
=E_{rr}.
\end{equation}
Thus $\Curv(p_i,u_{ij})=-q_ip_j$, $\Curv(p_i,v_{ir})=
-E_{rr}$.
Let $v_j$ be of the form $u_{jk}$ with $j\geq k$ 
or $v_{jr}$. Then
\[
\ev\Theta_{2n}(p_1\wedge v_1\wedge\cdots
\wedge p_n\wedge v_n)
=
(-1)^n
Q_n\left(\frac{\partial v_1}{\partial q_1},
\dots,\frac{\partial v_n}{\partial q_n}\right),
\]
as the only permutations contributing non-trivially are
those for which $\sigma(2)=\sigma(1)+1$, $\sigma(4)=
\sigma(3)+1$ and so on. By \eqref{e-PeterWalsh} this
identity gives a formula for $Q_{n}$ in terms of
$\Theta_{2n}$, as it determines the value of $Q_n$
on $\mathfrak h$ uniquely.
On the other hand, the left-hand side can be computed
using the definition of $\Theta_{2n}$:
by evaluating $\pi_{2n}$ on
$p_1\otimes v_1\otimes\cdots\otimes p_n\otimes v_n$,
only one permutation contributes and
we obtain (see \eqref{e-omega})
\begin{eqnarray*}
\lefteqn{
\ev\Theta_{2n}(p_1\wedge v_1\wedge
\dots\wedge p_n\wedge v_n)}\\
&=&
\mu_{2n}\circ\int_{\Delta_{2n}}\prod_{1\leq i\le j\leq n}
\omega_{2n}
\mathrm{tr}\,\left(
1\otimes 1\otimes \frac{\partial{v_1}}{\partial q_{1}}
\otimes\cdots\otimes 1\otimes 
\frac{\partial{v_n}}{\partial q_{n}}\right)
\pm\cdots,
\end{eqnarray*}
where the omitted terms are obtained by permuting the last $2n$ factors.
By Lemma \ref{l-103}, the effect of these omitted term is to replace
the integration over the simplex $\Delta_{2n}$ with the integration over
$[0,1]^{2n}$. The integrand depends then only on $n$ variables, so
the integral reduces to an integral over $[0,1]^n$. The claim follows.
\hfill $\square$

\medskip

We now compute $P_n$ on $\mathfrak d$ by evaluating the
integral in \eqref{e-Purvis}. 

To get the signs right, notice that
the Lie algebra $\mathfrak{gl}_n$ is embedded in $W_{n,N}$ via $(a_{ij})\mapsto
-\sum a_{ij}q_ip_j$ and $\mathfrak{gl}_N$ is embedded
as the subalgebra of constant matrices. Thus a
diagonal matrix $\mathrm{diag}(t_1,\dots,t_n)\oplus
\mathrm{diag}(s_1,\dots,s_N)\in \mathfrak{gl_n}\oplus\mathfrak{gl}_N$
corresponds to 
\[X=-\sum t_iq_ip_i+\sum s_r E_{rr}\in 
\mathfrak d.
\]
We compute $P_n(X,\dots,X)$.
Since $X$ is of degree
at most 2 in $p_i,q_i$ we may replace the product of exponentials
in the definition of $P_n$ by 
\begin{equation}\label{e-lR1}
\prod_{i<j}(1+
\epsilon\psi(u_i-u_j)\alpha_{ij}
+\frac{\epsilon^2}2\psi(u_i-u_j)^2\alpha_{ij}^2).
\end{equation}
 Expanding the product gives a
sum of terms labeled by graphs with vertex set $\{1,
\dots,n\}$, with an edge between $i$ and $j$ indicating a factor
$\psi(u_i-u_j)\epsilon\alpha_{ij}$ in the integrand. Only graphs with zero
or two edges emerging from each vertex
give a non-trivial contribution to $P_n$.  
Their connected components
are thus cycles. The integral over $u_1,\dots,u_n$ can then be written
as a product of integrals, one for each cycle: a cycle of length $j$
gives a contribution
\[
I_j=\int_{[0,1]^j}\psi(u_1-u_2)\psi(u_2-u_3)\cdots\psi(u_{j}-u_1)
du_1\cdots du_j.
\]
This vanishes for reason of symmetry if $j$ is odd.
Altogether we then have a sum labeled by the numbers $\ell_j\geq 0$ of 
connected components with $j$ vertices $(j\geq1)$:
\[
P_n(X,\dots,X)
=
\sum_{\sum_jj\ell_j=n}C_{\ell_1,\ell_2,\dots}
\prod_{j\geq 2}
\left(\frac{2I_j}{2^{j}}\sum_{i=1}^\infty 
(\epsilon t_i)^j\right)^{\ell_j}
\sum_{r=1}^N s_r^{\ell_1},
\]
where we set $t_i=0$, $i>n$.
The combinatorial factor is
\[
C_{\ell_1,\ell_2,\dots}=
\frac
{n!}
{{\ell_1!}\prod_{j\geq2}\ell_{j}!(2j)^{\ell_j}}.
\]
It is computed as follows. All graphs with $\ell_j$ cycles of length $j$
are obtained from any one of them by a permutation of vertices. The
denominator is $2^{\ell_2}$ (coming from the factor $1/2$ in \eqref{e-lR1})
times the order of the automorphism subgroup. The latter is generated
by  cyclic permutations and reflections
of vertices in each cycle (giving a factor $2j$ for each cycle of length
$j\geq2$ and $2$ for each factor of length 2)
 and the permutations of cycles of equal lengths (giving the factor
$\prod\ell_j!$).

We now consider the generating function
\[
S=\sum_{n=0}^\infty \frac1{n!}P_n(X,\dots, X),
\]
(we set $P_0=1$) as a formal power series in $t_1,t_2,\dots, s_1,\dots,s_N$.
We have
\[
S=\sum_{r=1}^N\exp\left(s_r+
\sum_{j\geq2}
\frac {I_j}{2^jj}\sum_{i=1}^\infty (\epsilon t_i)^j\right).
\]
\begin{lemma}
\[I_j=\left\{
\begin{array}{ll}
0,& \text{$j$ odd,}\\
2(\pi \ii)^{-j}\sum_{k=1}^\infty\frac1{k^{j}},& \text{$j$ even.}
\end{array}
\right.
\]
\end{lemma}

\noindent{\it Proof:}
The integral can be evaluated by Fourier series, using
\[
\psi(u)=\sum_{k\in\mathbb Z-\{0\}}\frac{\ii}{\pi k}e^{2\pi \ii k u}.
\]
\hfill$\square$

\medskip

By using this formula, it is easy to evaluate the sum over $j$ appearing in
the formula for $S$:
\begin{eqnarray*}
\sum_{j=2}^\infty\frac{I_j}{2^jj}t^j&=&-\sum_{k=1}^\infty\ln
\left(1-\left(\frac {\epsilon t}{2\pi \ii k}\right)^2\right)\\
&=&\ln \frac {\epsilon t/2}{\sinh(\epsilon t/2)}.
\end{eqnarray*}
Putting everything together, we finally get
\[
\sum_{n=0}^\infty \frac1{n!}P_n( X,\dots,X)
=
\prod_{i=1}^\infty\frac {\epsilon t_i/2}
{\sinh(\epsilon t_i/2)}\;\sum_{r=1}^Ne^{s_r}.
\]
Thus for general $X=X_1\oplus X_2\in\mathfrak{sp}_{2n}\oplus \mathfrak{gl}_N$,
\[
P_n(X,\dots,X)
=\left[\det \left(\frac{\epsilon X_1/2}{\sinh(\epsilon X_1/2)}
\right)^{\frac12}
 \mathrm{tr}\,e^{X_2}\right]_n.
\]
This completes the proof of Theorem \ref{t-lRRH}.
\hfill
$\square$

\subsection{The trace of $1$}\label{ss-Clarissa}
We return to the setting of Sect.~\ref{sect3} and compute the
trace of $1\in \mathcal A_D(M)$ for a compact symplectic manifold
$M$. By Theorem \ref{t-2}, we have
\[
\mathrm{Tr}_D(1)=\frac 1{(2\pi\ii\epsilon)^n}\int_M\psi_D(1)
=\frac {(-1)^n}{(2\pi\ii\epsilon)^n(2n)!}\int_M\ev\Theta_{2n}(A^{2n}).
\]
We now use Theorem \ref{t-lRRH}. Set $\mathfrak g=\mathcal A^\pol_{2n}$,
regarded as a Lie algebra with bracket $[\ ,\ ]_\epsilon$,
$\mathfrak{h}=\mathfrak{sp}_{2n}\oplus\mathfrak {gl}_1\subset \mathfrak g$.
Let $\pr:\mathfrak g\to
\mathfrak h$ be the $\mathfrak h$-equivariant projection
onto the quadratic and constant part:
\[
\pr(f(y))=
f(0)+\sum\frac12\partial_i\partial_jf(0)y_iy_j.
\]
Then the curvature $\Curv$ (see \ref{ss-MrsLucreziaWarrenSmith})
can be computed.
Since $A$ is defined up to a central 1-form and
the quadratic part may be absorbed in the connection,
we may assume that $\pr(A)=0$. Since $\nabla$
preserves the degree in $y$, we have $0=\nabla\pr(A)=\pr\nabla A$.
For any
vector fields $\xi,\eta$ on $M$, we  obtain
(see Section \ref{ss-Sally})
\begin{eqnarray*}
C(A(\xi),A(\eta))&=&-\pr([A(\xi),A(\eta)]_\epsilon)\\
&=&
-\pr(\nabla A(\xi,\eta)+[A(\xi),A(\eta)]_\epsilon)\\
&=&
\tilde R(\xi,\eta)-\Omega(\xi,\eta)
\end{eqnarray*}
Let $P_n=(\hat{\mathrm{A}}_\epsilon\,\mathrm{Ch})_n
\in (S^h\mathfrak h)^{*\mathfrak h}
$ be the
polynomial on the right-hand side of 
Theorem \ref{t-lRRH}, so that
$\ev\Theta_{2n}=(-1)^n\chi(P_n)$.
Then for any vector fields $\xi_1,\dots,\xi_{2n}$,
\begin{eqnarray*}
\lefteqn{
\frac{(-1)^n}{(2n)!}\ev\Theta_{2n}(A^{2n})
(\xi_1,\dots,\xi_{2n})}
\\
&&=
(-1)^n\ev\Theta_{2n}
(A(\xi_1)\wedge\cdots\wedge A(\xi_{2n}))
\\
&&=
\frac1{n!}\sum_{\sigma(2i-1)<\sigma(2i)}
\mathrm{sign}(\sigma)
P_n(C(\xi_{\sigma(1)},\xi_{\sigma(2)})\cdots C(\xi_{\sigma(2n-1)},\xi_{\sigma(2n)}))
\\
&&=
\frac1{n!}P_n((R-\Omega)^n)(\xi_1,\dots,\xi_{2n}).
\end{eqnarray*}
By homogeneity, $(\hat{\mathrm{A}}_\epsilon(R)\mathrm{Ch}(
-\Omega))_n=
\epsilon^n (\hat{\mathrm{A}}(R)\mathrm{Ch}(
-\epsilon^{-1}\Omega))_n$.
We obtain a formula for the trace of 1:
\thm (Fedosov, Nest--Tsygan) Let $R$ be the curvature 
of the symplectic connection $\nabla$. Then
\begin{eqnarray*}
\mathrm{Tr}_D(1)=\frac1{(2\pi\ii)^n}\int_M\hat{\mathrm A}(R)
\exp(-\Omega/\epsilon)
\end{eqnarray*}

\end{document}